\documentclass[oneside, 12pt]{amsart}
\usepackage{amscd, amssymb, amsmath, mathrsfs}
\usepackage[english]{babel}
\usepackage{booktabs}
\usepackage{tikz-cd}
\usepackage[all]{xy}
\usepackage[pdftex, colorlinks=true,  citecolor=blue, linkcolor=blue]{hyperref}

\newcommand\mylabel[1]{\label{#1}\marginpar{\vspace{-1ex}\medskip\medskip\footnotesize \tt #1}}
\renewcommand\mylabel[1]{\label{#1}}
\newcommand{\mydate}{
\number\day\space
\ifcase\month \or January\or February\or March\or April\or May\or June\or July\or August\or September\or October\or November\or December\fi 
\space\number\year}

\DeclareUrlCommand\arXiv{\urlstyle{same}}

\newtheorem{theorem}{Theorem}[section]
\newtheorem{maintheorem}{Theorem}
\newtheorem{lemma}[theorem]{Lemma}
\newtheorem{proposition}[theorem]{Proposition}
\newtheorem{corollary}[theorem]{Corollary}

\newtheorem*{conjecture*}{Conjecture}

\theoremstyle{definition}
\newtheorem{definition}[theorem]{Definition}

\newtheorem*{acknowledgement}{Acknowledgement}

\theoremstyle{remark}

\newcommand{\ZZ}{\mathbb{Z}}
\newcommand{\QQ}{\mathbb{Q}}

\newcommand{\CC}{\mathbb{C}}
\newcommand{\FF}{\mathbb{F}}

\newcommand{\PP}{\mathbb{P}}

\newcommand{\GG}{\mathbb{G}}

\newcommand{\shA}{\mathscr{A}}

\newcommand{\shE}{\mathscr{E}}
\newcommand{\shF}{\mathscr{F}}

\newcommand{\shN}{\mathscr{N}}
\newcommand{\shL}{\mathscr{L}}

\newcommand{\foA}{\mathfrak{A}}
\newcommand{\foB}{\mathfrak{B}}

\newcommand{\foX}{\mathfrak{X}}

\newcommand{\aff}{\text{\rm aff}}
\newcommand{\Aff}{\text{\rm Aff}}
\newcommand{\alg}{\text{\rm alg}}

\newcommand{\Alb}{\operatorname{Alb}}

\newcommand{\Aut}{\operatorname{Aut}}

\newcommand{\can}{\text{\rm  can}}

\newcommand{\Dom}{\operatorname{Dom}}

\newcommand{\End}{\operatorname{End}}

\newcommand{\Et}{{\text{\rm Et}}}

\newcommand{\Frac}{\operatorname{Frac}}

\newcommand{\Gal}{\operatorname{Gal}}
\newcommand{\GL}{\operatorname{GL}}

\newcommand{\Hom}{\operatorname{Hom}}

\newcommand{\id}{{\operatorname{id}}}
\newcommand{\Image}{\operatorname{Im}}

\newcommand{\Kernel}{\operatorname{Ker}}

\newcommand{\Lie}{\operatorname{Lie}}

\newcommand{\invlim}{\varprojlim}

\newcommand{\lra}{\longrightarrow}

\newcommand{\Mat}{\operatorname{Mat}}

\newcommand{\maxid}{\mathfrak{m}}

\renewcommand{\O}{\mathscr{O}}

\newcommand{\op}{\text{\rm op}}
\newcommand{\ord}{\operatorname{ord}}

\newcommand{\perf}{{\text{\rm perf}}}
\newcommand{\Pic}{\operatorname{Pic}}

\newcommand{\pr}{\operatorname{pr}}

\newcommand{\quadand}{\quad\text{and}\quad}

\newcommand{\ra}{\rightarrow}

\newcommand{\rank}{\operatorname{rank}}
\newcommand{\red}{{\operatorname{red}}}

\newcommand{\Res}{\operatorname{Res}}

\newcommand{\sep}{{\operatorname{sep}}}
\newcommand{\Set}{{\text{\rm Set}}}

\newcommand{\Sp}{\operatorname{Sp}}
\newcommand{\Spec}{\operatorname{Spec}}

\newcommand{\val}{\operatorname{val}}

\newcommand{\uHom}{\underline{\operatorname{Hom}}}

\begin{document}

\title[]
      {Varieties with free tangent sheaves}

\author[Damian R\"ossler]{Damian R\"ossler}
\address{Mathematical Institute, University of Oxford, Andrew Wiles Building, Rad-
cliffe Observatory Quarter, Woodstock Road, Oxford, OX2 6GG, United Kingdom}
\curraddr{}
\email{rossler@maths.ox.ac.uk}

\author[Stefan Schr\"oer]{Stefan Schr\"oer}
\address{Heinrich Heine University D\"usseldorf, Faculty of Mathematics and Natural Sciences, Mathematical Institute, 40204 D\"usseldorf, Germany}
\curraddr{}
\email{schroeer@math.uni-duesseldorf.de}

% 14K05 Algebraic theory of abelian varieties
% 14J50 Automorphisms of surfaces and higher-dimensional varieties
% 14J27 Elliptic surfaces, elliptic or Calabi-Yau fibrations
% 14L15 Group schemes
% 14G17 Positive characteristic ground fields in algebraic geometry

\dedicatory{14 March 2025}

\iffalse
J. Reine Angew. Math.
Daniel Huybrechts

Forum Math. Sigma
Vakil, Totaro

Math. Ann.
Vasudevan Srinivas
Ordinary varieties with trivial canonical bundle are not uniruled
A characterization of ordinary abelian varieties by the Frobenius push-forward of the structure sheaf.

Compositio:
Mustata
Projective varieties with nef tangent bundle in positive characteristic

Duke
JEMS
JAG
AG
\fi

\begin{abstract}
We coin the term \emph{$T$-trivial varieties} to denote smooth proper schemes over   ground fields $k$ whose tangent sheaf is free.
Over the complex numbers, this are precisely the abelian varieties.
However, Igusa observed that in characteristic $p\leq 3$ certain bielliptic surfaces are $T$-trivial.
We show that   $T$-trivial varieties $X$ separably dominated by   abelian varieties $A$ can exist only for $p\leq 3$.
Furthermore, we prove that every $T$-trivial variety, after passing to a finite \'etale covering,
is fibered in  $T$-trivial varieties with Betti number $b_1=0$. 
We also show that if some $n$-dimensional $T$-trivial $X$ lifts to characteristic zero and  $p\geq 2n+2$ holds,
it admits a finite \'etale covering by an abelian variety. Along the way, we   establish several
results about the automorphism group of abelian varieties, and the existence of relative Albanese maps.

\end{abstract}

\maketitle
\tableofcontents

\newcommand{\dashra}{\dashrightarrow}
\newcommand{\lieg}{\mathfrak{g}}
\newcommand{\lieh}{\mathfrak{h}}
\newcommand{\PS}{{\PP^n}}
\newcommand{\frob}{{(p)}}
\newcommand{\uExt}{\underline{\operatorname{Ext}}}
\newcommand{\liea}{\mathfrak{a}}

%===========================================================
\section*{Introduction}
\mylabel{Introduction}

Let $k$ be a ground field of characteristic $p\geq 0$, and $X$ be a smooth proper scheme
with $h^0(\O_X)=1$, of dimension $n\geq 0$. Recall that the term \emph{$K$-trivial variety } often refers
to those $X$ where  he dualizing sheaf $\omega_X=\det(\Omega^1_{X/k})$ is isomorphic to the structure sheaf.
In this paper we are interested in the much stronger condition, where the sheaf of K\"ahler differentials $\Omega^1_{X/k}$
or equivalently the tangent sheaf $\Theta_{X/k}=\uHom(\Omega^1_{X/k},\O_X)$ 
itself are isomorphic to $\O_X^{\oplus n}$, in other words, these locally free sheaves are \emph{free}.
We find it convenient to coin the term \emph{$T$-trivial varieties} for such $X$. 

To simplify exposition, we assume  that $k$ is algebraically closed in the following discussion.
For every abelian variety $A$, the group law gives an identification
of the tangent sheaf with   $\Lie(A)\otimes_k\O_A$, and consequently such schemes are  $T$-trivial.
In characteristic zero, the theory of Albanese maps easily shows that every $T$-trivial variety arises
in this way.  However, Igusa \cite{Igusa 1955} noted that  for $p\leq 3$   certain bielliptic surfaces  $X=(E\times E')/G$  become $T$-trivial.
Building on this, 
Li \cite{Li 2010} conjectured that  the $T$-trivial varieties in characteristic $p\geq 5$ are precisely the $T$-trivial varieties.
Our paper is concerned with the following natural questions:

\begin{enumerate}
\item Suppose a $T$-trivial variety $X$   that is not an abelian variety arises as a quotient $X=A/G$ of an abelian variety.
Does this imply $p\leq 3$?
\item What can be said about the Albanese map $X\ra\Alb_{X/k}$ for $T$-trivial varieties $X$ in characteristic $p>0$?
\item 
Suppose a $T$-trivial variety $X$ in characteristic $p>0$ admits a lifting to characteristic zero. Does this imply that $X$ is an abelian variety?
\end{enumerate}

In their important work \cite{Mehta; Srinivas 1987}, Mehta and Srinivas    answered  (iii)  affirmatively, under the assumption that $X$ is ordinary and projective.
This was refined by Joshi \cite{Joshi 2021}, who checked     in characteristic $p\geq 5$ that 
every $T$-trivial surface   is an  abelian surface,  and by Li \cite{Li 2010}, 
who established that for  $p\geq 3$ that every ordinary $T$-trivial variety is an  abelian variety.  

Roughly speaking, our  contributions   are as follows: First,  we establish that (i) indeed is true.
Second, we show that every $T$-trivial variety $X$ admits a finite \'etale covering $\tilde{X}$ where the fibers of the  Albanese map  
are $T$-trivial varieties with Betti number $b_1=0$. Third, we show that (iii) holds if characteristic and dimension satisfy   $p\geq 2\dim(X)+2$.

Let us now describe these results in more detail. Our first main results clarifies if and in  which ways Igusa's
construction $X=(E\times E')/G$ might be carried out at other primes  or in higher dimensions.

\begin{maintheorem} 
(See Thm.\ \ref{t-trivial with etale covering})
Let $X$ be a $T$-trivial variety that is not an abelian variety, but has  a finite   surjection $A\ra X$ 
from an abelian variety, with   $k(X)\subset k(A)$ separable.  Then the following holds:
\begin{enumerate}
\item The characteristic satisfies $p\leq 3$. 
\item The  abelian variety $A$ is not   simple and contains a   point of order $p$.
\item For $p=3$ the $T$-trivial variety $X$ is not ordinary, and the abelian variety $A$ has a supersingular quotient.  
\end{enumerate}
\end{maintheorem}

Indeed,  for Igusa's bielliptic surface $X=(E\times E')/G$  one  needs on $A=E\times E'$ a point of order $p$,
and for $p=3$ it turns out that the factor $E$ must be supersingular.
Our  result also answers a question of Joshi (\cite{Joshi 2021}, Section 6).
%\mynote{which one?}
Furthermore, one obtains a  new proof for Li's Theorem mentioned above (\cite{Li 2010}, Theorem 4.2).

Note that quotients $X=A/G$ as above are sometimes called \emph{hyperelliptic varieties} \cite{Lange 2001}.
Over the complex numbers, their classification problem  boils down to understand discontinuous groups
of affine transformations inside the semidirect product $\CC^n\rtimes\GL_n(\CC)$, a topic that received considerable attention  in dimension three
(\cite{Uchida; Yoshihara 1976}, 
\cite{Catanese; Demleitner 2018}, \cite{Catanese; Demleitner 2020}, \cite{Catanese; Demleitner 2023}, \cite{Bauer; Gleissner 2022}).

The proof for the above theorem relies  on a statement on   automorphism groups of abelian varieties, which seems
to be of independent interest:

\begin{maintheorem} 
(See Thm.\ \ref{p-group})
Let $A$ be an abelian variety in characteristic $p>0$. Suppose
the kernel of $\Aut(A)\ra\GL(\Lie(A))$ contains a non-trivial element $h$ of finite order.
Then the following holds:
\begin{enumerate}
\item The characteristic must be $p\leq 3$.
\item The order of $h\in\Aut(A)$ is a $p$-power.
\item For $p=3$ the abelian variety $A/\Kernel(\id-h)$ is supersingular.
\end{enumerate}
\end{maintheorem}

This can be seen as a variant of Minkowski's Theorem (\cite{Minkowski 1887}, Section 1),
which states that an integral matrix  $H\in\Mat_n(\ZZ)$   with $H\equiv E$ modulo 
some integer $r\geq 3$ is already the identity matrix, and Serre's result (\cite{Grothendieck 1961}, Appendix), which asserts that 
an automorphism $h\in\Aut(A)$ that is the identity on the group scheme $A[r]$ for some $r\geq 3$ is already the identity.
Our proof relies  on the \emph{Weil Conjectures} for abelian varieties over prime fields, and   number-theoretical properties
of the eigenvalues of Frobenius, the so-called \emph{Weil numbers}.
 
Our third main results implies that if Li's Conjecture fails, it must fail in a spectacular way:

\begin{maintheorem} 
(See Thm.\ \ref{decomposition t-trivial})
For every  $T$-trivial variety $X$, there is a finite \'etale covering $X'\ra X$ such that the fibers
of the Albanese map $X'\ra\Alb_{X'/k'}$ are $T$-trivial varieties with Betti number $b_1=0$.
\end{maintheorem}

The $T$-trivial varieties  with $b_1=0$ would be extremely remote from   the situation 
over the complex numbers,   and from  any  Igusa-type construction in positive characteristics as well.
Except for singletons, we have no clue so far whether or not such bizarre objects exist.
 
For the proof for the above result we use the theory of  \emph{relative Albanese map}.
It is now high  time  to dismiss the assumption that the ground field $k$ is algebraically closed.
Recall that a smooth proper scheme $P$  with $h^0(\O_P)=1$ is called \emph{para-abelian variety} if, for some 
field extension $k\subset k'$, the base-change $P'=P\otimes k'$ can be endowed with  a group law, and thus becomes
an abelian variety. This notion was developed in \cite{Laurent; Schroeer 2024} and \cite{Schroeer 2024}, and actually goes back to Grothendieck \cite{FGA VI}.
It turns out that $G=\Aut^0_{P/k}$ is an abelian variety, acting freely and transitively on $P$. 
Summing up, para-abelian varieties allow  for an  intrinsic way to handle torsors with respect to abelian varieties.
In connection to $T$-trivial varieties, which in the first place have no distinguished point, it is indeed preferable to work with para-abelian varieties
instead of abelian varieties.
 
For proper flat morphism $f:X\ra S$ of finite presentation, the \emph{relative Albanese variety}
is a family of para-abelian varieties $\Alb_{X/S}$,
and the \emph{relative Albanese map} is a universal arrow $X\ra \Alb_{X/S}$ to such families.
Their existence depends on particular properties of $\Pic_{X/S}$, which do not always hold (\cite{Laurent; Schroeer 2024}, Theorem 10.2).
In Corollary \ref{unconditional existence}, we   provide a new   unconditional statement in characteristic zero. 
For   $T$-trivial varieties $X$ with $p>0$, we seek to form the  relative Albanese variety   with respect to the
absolute Albanese map, but in this setting existence is unclear. The following work-around, which relies
on the \emph{Weil  Extension Theorem},   seems to be of independent interest:

\begin{maintheorem} 
(See Thm.\ \ref{weak relative albanese map})
Suppose that  $S$ is normal, and that the generic fiber $X_\eta$ contains a rational point. 
After removing a closed set $Z\subset S$ of codimension at least two,
$P_\eta=\Alb_{X_\eta/\kappa(\eta)}$ extends
to a family of abelian varieties $P$ over $S$, and the Albanese map $ g_\eta:X_\eta\ra P_\eta$ extends to a morphism $g:X\ra P$.
\end{maintheorem}

Our last main result connects the theory of $T$-trivial varieties with lifting properties.
Suppose $X$ is smooth and projective over an algebraically closed ground field  of characteristic $p\geq 0$.
For $k=\CC$ it follows from Yau's proof of the Calabi Conjecture that $X$ admits a finite \'etale covering by an abelian variety if and only 
the Chern classes $c_1$ and $c_2$ vanish. It is unclear to what extend the reverse implication holds true for $p>0$.
We show:

\begin{maintheorem} 
(See Thm.\ \ref{liftings and coverings}) 
In the above situation, suppose the following holds:
\begin{enumerate}
\item For some $\ell\neq p$, the  $\ell$-adic Chern classes $c_1$ and $c_2$   both vanish.
\item The scheme $X$ projectively lifts to characteristic zero.
\item Characteristic and dimension satisfy  $p\geq 2n+2$.
\end{enumerate}
Then there is a finite \'etale covering $A\ra X$ by some  abelian variety $A$. 
\end{maintheorem}

Note that condition (i) automatically  holds if $X$ is $T$-trivial.
In dimension $n=2$, we see that  every $T$-trivial surface in characteristic $p\geq 7$ arises from an abelian surface. Note that
this already actually holds  for $p\geq 5$, by the  
Bombieri--Mumford generalization of the Bagnera--de Francis classification for bielliptic surface
(\cite{Bombieri; Mumford 1977}, Section 3).

The paper is organized as follows: In Section \ref{Generalities} we collect   generalities on Weil restriction, free sheaves, abelian varieties,
$\ell$-adic cohomology, and algebraic fundamental groups. Section \ref{Automorphisms trivial} contains results on 
automorphisms of abelian varieties that act trivially on the Lie algebra.
We introduce the  notion of $T$-trivial varieties $X$ and establish their basic properties in Section \ref{T-trivial}.
There we also establish some structure results if it is dominated by some abelian variety.
In Section \ref{Relative albanese} the theory of relative Albanese maps is developed further.
This is used in Section \ref{Betti number} to obtain a splitting result where $T$-trivial varieties with $b_1=0$ appear.
In the final Section \ref{Liftability}, we connect the theory of $T$-trivial varieties with liftings to characteristic zero.

\begin{acknowledgement}
Both authors are grateful for the hospitality during their mutual  visits at the Heinrich Heine University D\"usseldorf
and    Pembroke College Oxford.
The research was supported by the Deutsche Forschungsgemeinschaft via the grant  \emph{SCHR 671/10-1 Varieties with Free Tangent Sheaf}. 
It was also conducted       in the framework of the   research training group
\emph{GRK 2240: Algebro-Geometric Methods in Algebra, Arithmetic and Topology}.
\end{acknowledgement}

%===========================================================
\section{Generalities}
\mylabel{Generalities}

In this section we collect some general facts on Weil restriction,   free sheaves,
abelian varieties,   $\ell$-adic cohomology and algebraic fundamental groups that are relevant throughout, and perhaps
of independent interest. For simplicity, we work over a ground field $k$, of characteristic $p\geq 0$.

Let us start with \emph{Weil restrictions}. Suppose 
$k_0\subset k$ is a subfield such that the   degree of the extension is finite.
For each  scheme or algebraic space  $X$ over $k$, the functor 
$$
(\Aff/k_0)\lra(\Set),\quad R_0\longmapsto X(R_0\otimes_{k_0}k)
$$
is called   \emph{Weil restriction} $X_0=\Res_{k/k_0}(X)$. 
It is representable by an algebraic space  (\cite{Ji; Li; McFaddin; Moore; Stevenson 2022}, Theorem 6.5.2), which we denote by the same symbol.
Note that $X_0$ is schematic provided that $X$ has the \emph{AF-property}, which means that
every finite set of points admits a common affine neighborhood. In any case,   Weil restriction is right adjoint 
to base-change, such that    
$\Hom(Y,\Res_{k/k_0}(X))=\Hom(Y\otimes_{k_0} k,X)$.
The following property is most useful (loc.\ cit., Theorem 6.1.5):

\begin{lemma}
\mylabel{weil restriction}
Suppose   $k_0\subset k$ is separable, with Galois closure  $k'$. Then 
\begin{equation}
\label{pullback weil restriction}
\Res_{k/k_0}(X)\otimes_{k_0} k' = \prod_{\iota:k\ra k'}  (X\otimes_kk'),
\end{equation} 
where the product runs over all $k_0$-embeddings $\iota:k\ra k'$. 
\end{lemma}

If $k_0\subset k$ is already Galois, the indices $\iota$ become the  elements from the Galois group.
In any case, the $\iota$ form a set of cardinality $[k:k_0]$, according to 
\cite{A 4-7}, Chapter V, \S6, No.\ 3, Corollary to Proposition 1,  giving the dimension formula
$$
\dim(X_0)=[k:k_0]\cdot \dim(X).
$$
Moreover,  $X_0=\Res_{k/k_0}(X) $ is  the quotient of  \eqref{pullback weil restriction} by the canonical  action of $G'=\Gal(k'/k_0)$, which is free.
Furthermore, we see by descent that if $X$ is   smooth and proper   with $h^0(\O_X)=1$, the same holds for the Weil restriction
$X_0=\Res_{k/k_0}(X)$.
Also note that the identification \eqref{pullback weil restriction} is natural: 
Given $f:X\ra Y$, the induced morphism $f_0=\Res_{k/k_0}(f)$ on Weil restrictions
has the property 
\begin{equation}
\label{pullback natural}
f_0\otimes_{k_0} k'=\prod_{\iota:k\ra k'} (f\otimes_k k') 
\end{equation} 

We next turn to    \emph{free sheaves}. 
Recall that a  locally free sheaf $\shE$ of rank $r\geq 0$ on  a scheme  $X$ is called \emph{free}
if it is isomorphic to $ \bigoplus_{i=1}^ r\O_X  $. 
The following categorical observation is   useful:

\begin{lemma}
\mylabel{categorical view}
Suppose $h^0(\O_X)=1$. Then the functor $V\mapsto V\otimes_k\O_X$ from the abelian category finite-dimensional $k$-vector spaces $V$ to the 
abelian category of quasicoherent sheaves $\shE$ on $X$ is exact and fully faithful, and its essential image
is the category of free sheaves of finite rank.
\end{lemma}

\proof
The statement on the essential image is obvious.  
The functor is fully faithful, because both $\Hom(k^n,k^m) $ and $\Hom(\O_X^{\oplus n},\O_X^{\oplus m})$
are given by $m\times n$-matrices $(\varphi_{ij})$ with entries from $\Hom(k,k)=k=\Gamma(X,\O_X)=\Hom(\O_X,\O_X)$.
The functor is exact since the structure map $X\ra \Spec(k)$ is flat. 
\qed

\medskip
It follows that  all short exact sequences $0\ra \shE'\ra \shE\ra\shE''\ra 0$ with free sheaves of finite rank  are  split.
Moreover, a map  $\varphi:\shE\ra\shF$ of free sheaves of finite rank is injective  or surjective provided that corresponding property
holds after tensoring with the residue field of  some   point $a\in X$.
The following useful observation is essentially contained in \cite{Novakovic 2012}, Lemma 4.2:

\begin{proposition}
\mylabel{descent freeness}
Suppose $h^0(\O_X)=1$, and let $\shE$ be quasicoherent sheaf on $X$. If for some field extension $k\subset k'$,
the base-change $\shE'=\shE\otimes k'$ to $X'=X\otimes k'$ becomes free of rank $r\geq 0$,
the same already holds for $\shE$.
\end{proposition}

\proof
We have $h^0(\shE) =  \dim_{k'}H^0(X',\shE')=\dim_{k'} H^0(X',\O_{X'})\cdot r=r$.
Choose a basis $s_1,\ldots,s_r\in H^0(X,\shE)$, and consider the resulting homomorphism of quasicoherent sheaves   
$$
k^r\otimes_k\O_X\stackrel{s_1,\ldots,s_r}{\lra} H^0(X,\shE)\otimes_k\O_X\stackrel{\can}{\lra} \shE.
$$
The map on the left is bijective by construction, and the map on the right becomes bijective after base-change to $k'$.
It follows that both maps are bijective, hence $\shE$ is free of rank $r$.
\qed

\medskip
Let us also record:

\medskip
\begin{proposition}
\mylabel{freeness after pullback}
Suppose $h^0(\O_X)=1$. Let  $\shE$ be locally free sheaf of finite rank  on $X$ that is globally generated,
and assume that there is a scheme $X'$ with $h^0(\O_{X'})=1$, and a surjection $f:X'\ra X$ such that the pullback $\shE'=f^*(\shE)$ is free.
Then $\shE$ is free.
\end{proposition}

\proof
Set $r=\rank(\shE)$ and $V=H^0(X,\shE)$. The canonical map $ V\otimes_k\O_X\ra\shE$ is surjective,
and the same holds for its pullback to $X'$.
By Lemma \ref{categorical view}, there is a  vector subspace $U\subset V$ such that the pullback of $\varphi:U\otimes_k\O_X\ra \shE$ is bijective.
We may regard   $s=\det(\varphi)$ as a section for the invertible sheaf $\shL=\det(\shE)$, and the task is to verify
that is has no zeros. Let $a\in X$ be a point and  $a'\in X'$ be some preimage. 
By construction   $s(a')=s(a)\otimes\kappa(a')$ does not vanish, hence $s(a)\neq 0$.
\qed

\medskip 
We next turn to \emph{abelian varieties} $A$ in positive characteristics $p>0$. 
Then $\Lie(A)$ is a  \emph{restricted Lie algebra}, having zero brackets $[x,y]=0$  and some \emph{$p$-map} $x\mapsto x^{[p]}$.
One says that $A$ is \emph{superspecial} if $\Lie(A)$ is isomorphic to $k^g$, the restricted Lie algebra
where both bracket and $p$-map are zero. In dimension $g=1$, this are precisely the supersingular elliptic curves.
Also note that up to twists, there is but one superspecial abelian variety
in each dimension $g\geq 2$, namely the product of supersingular elliptic curves,
a result attributed to Deligne (\cite{Shioda 1978}, Theorem 3.5, see also  \cite{Ogus 1979}, Theorem 6.2).
One calls $A$ \emph{supersingular} if and only if it is isogeneous to such a product $E_1\times\ldots\times E_g$
of supersingular elliptic curves. This condition can   be rephrased in terms of Dieudonn\'e modules,
and holds over $k$ if and only if it holds over $k^\alg$ 
(\cite{Oort 1974}, Theorem 4.2   and \cite{Yu 2020}, Theorem 1.2).
The following observation will be useful:

\begin{lemma}
\mylabel{supersingular quotient}
For the $g$-dimensional abelian variety   $A$,  the following are equivalent:
\begin{enumerate}
\item There is a non-zero supersingular quotient $A/N$.
\item Some abelian subvariety in  $A\otimes k^\alg$ has a non-zero supersingular quotient.
\end{enumerate}
Under these equivalent conditions,  the    map $H^g(A^{(p)},\O_{A^{(p)}})\ra H^g(A,\O_A)$ induced by the relative Frobenius
$F:A\ra A^{(p)}$ is zero.
\end{lemma}

\proof
The implications (i)$\Rightarrow$(ii) is  trivial.
For the converse,  let $k\subset k'$ be any field extension, set $A'=A\otimes k'$, and consider the ordered set of abelian subvarieties
$N'_\lambda\subset A'$, $\lambda\in L$ such that $A'/N'_\lambda$ is supersingular.
For any two members $N'_\lambda$ and $N'_\mu$, the abelian variety $(N'_\lambda\cap N'_\mu)^0_\red$ belongs
to the family, because the class of supersingular abelian varieties are stable under products and subvarieties. 
In turn, our collection contains a smallest member $N'_0\subset A'$.
Obviously, this is stable under the action of $G=\Aut(k'/k)$. 

For $k'=k^\alg$, the Poincar\'e Irreducibility Theorem easily implies  that $A'/N'_0$ is non-zero.
Via Galois descent we see that  the base-change   to the perfect closure $k^\perf$
admits a non-zero supersingular quotient. 

Changing notation, we find a finite purely inseparable field extension
$k\subset k'$ with a non-zero supersingular quotient $f':A\otimes k'\ra B'$. The corresponding homomorphism $f:A\ra\Res_{k'/k}(B')$ 
factors over the abelian variety $\bar{A}=\Image(f)$. Pulling back we obtain
$$
A\otimes k'\lra \bar{A}\otimes k'\lra \Res_{k'/k}(B)\otimes k'\lra B.
$$
The map from $\bar{A}\otimes k'$ to $B$ is surjective, because $A\otimes k'\ra B$ is surjective, whereas the map to $\Res_{k'/k}(B)\otimes k'$
is a closed embedding, because this holds for $\bar{A}\ra\Res_{k'/k}(B)$.
Using that the kernel for the adjunction $\Res_{k'/k}(B)\otimes k'\ra B$ is affine (\cite{Conrad; Gabber; Prasad 2010}, Appendix A, Proposition 5.11), 
we conclude that $\bar{A}\otimes k'\ra B$ is an isogeny, and infer that $\bar{A}$ is supersingular.

For the remaining statement, suppose that
$B=A/N$ is a supersingular quotient. Passing to a further quotient,
we may assume that $B$ is supersingular elliptic curve, so the canonical map
$H^1(B^{(p)},\O_{B^{(p)}})\ra H^1(B,\O_B)$ induced by the relative Frobenius, which is identical to the Verschiebung
on $\Pic^0_{B/k}=B$, vanishes.
The cup product $\Lambda^g H^1(A,\O_A)\ra H^g(A,\O_A)$ is bijective (see for example
\cite{Roessler; Schroeer 2022}, Proposition 2.3).
Passing to the Stein factorization for the projection $A\ra B$, we may furthermore assume
that the canonical map $H^1(B,\O_B)\ra H^1(A,\O_A)$ is injective, and the same for Frobenius pullbacks.
Choose a non-zero $\epsilon_1\in H^1(B^{(p)},\O_{B^{(p)}})$ and
extend it to a basis $\epsilon_1,\ldots,\epsilon_g\in H^1(A^{(p)},\O_{A^{(p)}})$.
The relative Frobenius vanishes on the generator $\epsilon_1\cup\ldots\cup\epsilon_g$ of $H^g(A^{(p)},\O_{A^{(p)}})$,
because it vanishes the first factor.
\qed

\medskip
If the equivalent conditions of Lemma \ref{supersingular quotient} hold, we say that   $A$ \emph{has a supersingular quotient}.
This property played a crucial role in \cite{Pink; Roessler 2004}, where the term \emph{has a supersingular factor} was used.

Over finite ground fields  $k=\FF_q$, supersingularity can be characterized in terms of Frobenius eigenvalues:
Write $q=p^\nu$, and let $\Phi:A\ra A$ be the $\nu$-th power of the absolute Frobenius map.
Fix a prime $\ell\neq p$ and consider the induced $\QQ_\ell$-linear endomorphisms
on $H^i(A\otimes k^\alg,\QQ_\ell)$. For each embedding $\QQ_\ell\subset\CC$ the resulting eigenvalues
$\alpha_1,\ldots,\alpha_{2g}\in\CC$ are algebraic integers  and have absolute value $|\alpha_j|=p^{i/2}$,
according to \cite{Deligne 1980}, Corollary 3.3.9.
An algebraic integer  $\alpha\in\CC$ all whose conjugates $\alpha'$ have $|\alpha'|=p^{i/2}$ are called
\emph{Weil numbers of weight $i$}. Those of the particular simple
form $\alpha=\zeta\cdot p^{i/2}$ for some root of unity $\zeta$ are called \emph{supersingular}.
By \cite{Yu 2013}, Theorem 2.9 the abelian variety  $A$ is a supersingular  if and only if
the Weil numbers  $\alpha_1,\ldots,\alpha_{2g}\in\CC$  are supersingular.

We next come to  \emph{$\ell$-adic cohomology}, where $\ell>0$ is a fixed prime that is  invertible in the ground field $k$.
Let $X$ be  a   scheme, for the sake of exposition assumed to be proper.
Recall that  $\mu_{\ell^\nu}=\mu_{X,\ell^\nu}$, $\nu\geq 0$ denotes  the sheaf of $\ell^\nu$-th roots of unity
on the site $(\Et/X)$ of \'etale $X$-schemes, endowed with the \'etale topology.
The \emph{\'etale cohomology groups}
$H^i(X, \mu_{\ell^\nu}^{\otimes j})$, $\nu\geq 0$ form an inverse system of $\ZZ_\ell$-modules, and one defines the
\emph{$\ell$-adic cohomology} as 
$$
H^i(X,\ZZ_\ell(j)) = \invlim_\nu H^i(X,\mu_{\ell^\nu}^{\otimes j})\quadand H^i(X,\QQ_\ell(j)) = H^i(X,\ZZ_\ell(j)) \otimes_{\ZZ_\ell}\QQ_\ell.
$$
To discard    arithmetical contributions, one frequently considers the  cohomology groups  $H^i(X\otimes k^\alg,\ZZ_\ell(j))$  
and $H^i(X\otimes k^\alg, \QQ_\ell(j))$, which are finitely generated. The resulting \emph{Betti numbers} are   
$$
b_i(X) = \dim_{\QQ_\ell} H^i(X\otimes k^\alg, \QQ_\ell(j))
$$
Here both  $\ell$ and  $j$ are irrelevant, the latter because 
$H^i(X\otimes k^\alg,\ZZ_\ell(j))$ is obtained from $ H^i(X\otimes k^\alg,\ZZ_\ell)$ by tensoring with 
the invertible $\ZZ_\ell$-module $\ZZ_\ell(j)=\invlim_\nu (\mu_{\ell^\nu}(k^\alg)^{\otimes j})$. 
Also note that in all this one might use $k^\sep$ instead of $k^\alg$.

Although not via   cycle class maps, cohomology in degree one has the following well-known significance:

\begin{lemma}
\mylabel{betti and picard}
Suppose  $X$ is geometrically normal with $h^0(\O_X)=1$. Then the first Betti number is given by 
$b_1(X) = 2\dim(\Pic_{X/k})$.
\end{lemma}

\proof
First note  that  the Kummer sequence $0\ra \mu_{\ell^\nu}\ra\GG_m\stackrel{\ell^\nu}{\ra} \GG_m\ra 0$ yields 
an exact sequence
\begin{equation}
\label{from kummer}
0\lra k^\times/k^{\times \ell^\nu} \lra H^1(X,\ZZ/\ell^\nu\ZZ(1))\lra \Pic(X)[\ell^\nu]\lra 0.
\end{equation}
Next recall that the connected component $P=\Pic^\tau_{X/k}$ is a group scheme of finite type. It is actually proper (see for example \cite{Laurent; Schroeer 2024}, Proposition 2.3)
hence an extension of a finite group scheme $G$ by an abelian variety $A$. We thus have
 exact sequences $0\ra A[\ell^\nu]\ra P[\ell^\nu]\ra G[\ell^\nu]$ of group schemes, and obtain   exact sequences
$$
0\lra \invlim_\nu A[\ell^\nu](k) \lra \invlim_\nu P[\ell^\nu](k)\lra \invlim_\nu G[\ell^\nu](k) 
$$
of groups, where the terms on the right are finite. To proceed, we may assume that $k$ is algebraically closed.
Then the term on the left is a free $\ZZ_\ell$ module of rank $2\dim(A)$. The term in the middle 
can be identified with $H^i(X, \ZZ_\ell(1))$ by \eqref{from kummer}, and the assertion follows.
\qed

\medskip
Suppose now that $X$ is connected.
Let $a:\Spec(\Omega)\ra X$ be a geometric point, and $\pi_1(X,a)$ be the ensuing
\emph{algebraic fundamental group},    introduced by Grothendieck
(\cite{SGA 1}, Expos\'e V, Section 7) as     automorphism group of the fiber functor $V\mapsto V(\Omega)$,
which sends a finite \'etale covering $V$ to its fiber with respect to $a$. 
This functor indeed yields an equivalence between the category of finite \'etale coverings $V\ra X$
and the category of finite sets endowed with a continuous $\pi_1(X,a)$-action.

An abelian sheaf $F_\nu$ on the \'etale site $(\Et/X)$ is called
\emph{$\ell^\nu$-local system} if it is a twisted form of the constant sheaf
$(\ZZ/\ell^\nu\ZZ)^{\oplus r}_X$, for some rank $r\geq 0$. An \emph{$\ell$-adic local system} is an inverse system $F=(F_\nu)$,
where the entries are $\ell^\nu$-local systems, and the transition maps yield identifications
$F_\nu=F_\mu\otimes(\ZZ/\ell^\nu\ZZ)_X$ whenever $\mu\geq \nu$. 
Our geometric point $a:\Spec(\Omega)\ra X$ yields  the \emph{monodromy representation}
$\pi_1(X,a)\ra \GL( F_\nu\mid \Omega)$.
This actually gives an equivalence between the additive  category of $\ell^\nu$-local systems
of rank $r\geq 0$ and the additive category of  continuous representations
$\pi_1(X,a) \ra \GL_r(\ZZ/\ell^\nu\ZZ)$. For $\ell$-adic local systems $F=(F_\nu)$,
one obtains an equivalence to the additive category of continuous representations
$$
\pi_1(X,a)\lra \invlim_\nu\GL_r(\ZZ/\ell^\nu\ZZ) = \GL_r(\ZZ_\ell).
$$
Localizing the category of $\ell$-adic sheaves by tensoring the Hom sets with $\QQ_\ell$,
one obtains an equivalence to the category of continuous representations in $\GL_r(\QQ_\ell)$.
By abuse of notation, we also writes $\QQ_{\ell,X}$ for the $\ell$-adic local system $(\ZZ/\ell^\nu\ZZ)_X$, $\nu\geq 0$
in the localized category.

For each smooth proper morphism $f:X\ra Y$  and each $\ell^\nu$-local system $F_\nu$ on $X$,
the higher direct images $R^if_*(F_\nu)$ are $\ell^\nu$-local systems on $Y$. Moreover, 
for each cartesian diagram
$$
\begin{CD}
X'	@>>> 	X\\
@Vf'VV		@VVfV\\
Y'	@>>>	Y
\end{CD}
$$
the Proper Base-Change Theorem  yields an identification $R^if_*(F_\nu) | Y'  = R^if'_*(F_\nu|X')$,
with likewise  statements  for $\ell$-adic local systems $F=(F_\nu)$. 
 
Finally, suppose that $X$ is proper over the  ground field $k$, with $h^0(\O_X)=1$. Choose an algebraic closure $k^\alg$,
fix a closed point $x_0$ on the base-change $X\otimes k^\alg$, and set 
$S=\Spec(k)$. The structure morphism $X\ra S$ induces a short exact sequence
\begin{equation}
\label{homotopy sequence}
1\lra \pi_1(X\otimes k^\alg, x_0)\lra \pi_1(X,x_0)\lra \pi_1(S,x_0)\lra 1,
\end{equation} 
where the  term on the right becomes the Galois group $\Gal(k^\sep/k)=\Aut(k^\alg/k)$.
To simplify notation we set
$$
\Pi^\alg=\pi_1(X\otimes k^\alg, x_0)\quadand \Pi=\pi_1(X\otimes k^\alg, x_0)\quadand\Gamma=\pi_1(S,x_0).
$$
Conjugacy defines a homomorphism $\Gamma\ra\operatorname{Out}(\Pi^\alg)$. As explained in \cite{Brown 1982}, Chapter IV, Section 6,
the isomorphism classes of extensions with such an outer representation
become a principal homogeneous space with respect to the cohomology group $H^2(\Gamma, Z(\Pi^\alg))$, formed with respect to the center.

Recall that the open subgroups $H\subset\Pi$ are precisely
the closed subgroups of finite index. By Galois theory,   the transitive $\Pi$-set
$\Pi/H$  for such subgroups correspond to the finite \'etale covering $X'\ra X$ with connected total space. Moreover,
$$
\Gamma/H \quadand \bigcup_{\Pi^\alg\cdot\sigma\cdot H} \Pi^\alg/(\Pi^\alg\cap\sigma H\sigma^{-1})
$$
correspond to the \'etale $k$-algebra $k'=H^0(X',\O_{X'})$ and the base-change $X'\otimes k^\alg$,
respectively. The kernel $N\subset\Pi$ for the permutation representation on $\Pi/H$
is the largest normal subgroup contained in $H$, which is also open, and 
gives the Galois closure $X''$ for $X'$. Set $k'=H^0(X',\O_{X'})$ and $k''=H^0(X'',\O_{X''})$.
Note that we may easily have $h^0(\O_{X''})>h^0(\O_{X'})$, compare the discussion \cite{Debes; Douai 1997}, Section 2.8.
To get rid of such constant field extensions one may use the following observation:

\begin{lemma}
\mylabel{removing constant extensions}
The canonical morphism  $X''\ra X\otimes_kk''$    is a finite \'etale Galois covering, the $k''$-vector spaces
$H^0(X'',\O_{X''})$ and $H^0(X\otimes k'',\O_{X\otimes k''})$ are one-dimensional
\end{lemma}

\proof
Both schemes are \'etale over $X$, hence the morphism is \'etale by   \cite{EGA IVd}, Proposition 17.3.4.
Let $N\subset \Pi$ be the open normal subgroup corresponding to the composite map $X''\ra X$, and $\Gamma''\subset \Gamma$ be its image.
Then $X\otimes k''$ corresponds to the subgroup $\Pi''=\Pi\times_\Gamma\Gamma''$ of $\Pi$, and $X''\ra X\otimes_kk''$ corresponds to $N\subset\Pi''$.
The first statement follows.
The statement on $H^0(X'',\O_{X''})$ is trivial, and the statement on $H^0(X\otimes k'',\O_{X\otimes k''})$ follows from $h^0(\O_X)=1$.
\qed

%===========================================================
\section{Automorphisms of abelian varieties}
\mylabel{Automorphisms trivial}

The goal of this section is to establish several results  on automorphisms of abelian varieties, which will
play a crucial role in the next section, and  appears to be of independent interest.
 
Let $k$ be a ground field of characteristic $p\geq 0$, and $A$ be an abelian variety of dimension $g\geq 0$.
It comes with an associative algebra $\End(A)$, and  its unit group 
$\Aut(A)$ is a  countable  group. Each element  fixes the neutral element $e\in A$, and thus  stabilizes all
infinitesimal neighborhood $\Spec(\O_{A,e}/\maxid^{n+1})$. According to \cite{Matsumura 1963}, Lemma in Section 3 the resulting
linear representations
$$
\Aut(A)\lra \GL(\O_{A,e}/\maxid^{n+1})
$$
are injective for $n\geq 0$ sufficiently large.   Moreover, the kernels
of the above maps define a series of normal subgroups.
Throughout, we are interested in the case $n=1$, where the above   can also be seen as the canonical map
$\Aut(A)\ra \GL (\Lie(A))$. Our first main result reveals that the  torsion inside the kernel is rather restricted:

\begin{theorem}
\mylabel{p-group} 
Suppose $p>0$, and that the kernel of $\Aut(A)\ra\GL(\Lie(A))$ contains a non-trivial element $h$ of finite order.
Then the following holds:
\begin{enumerate}
\item The characteristic must be $p\leq 3$.
\item The order of $h\in\Aut(A)$ is a $p$-power.
\item For  $p=3$ the  abelian variety $A/\Kernel(h-\id)$ must be supersingular.
\end{enumerate}
\end{theorem}

\proof
The key idea is to understand the case where the ground field is a prime field, and then reduce to this situation via standard arguments.
We proceed in three steps.

\medskip\noindent
\textbf{Step 1:} \emph{Suppose that  $k=\FF_p$ is the prime field.}
First note that $A^{(p)}=A$ and that the relative and absolute Frobenius maps for $A$ coincide.
We thus have a short exact sequence $0\ra A[F]\ra A\stackrel{F}{\ra }A\ra 0$.
Also note that this sequence is functorial. Thus $F$ belongs to the center of the associative algebra
$\End(A)$, and every endomorphism stabilizes the Frobenius kernel $A[F]$.

Replacing our group element  $h\in\Aut(A)$ of finite order by a suitable power,
we may assume that $r=\ord(h)$ is prime.
The difference $ \id-h$ induces the zero map on $\Lie(A)$.
In light of the Demazure--Gabriel Correspondence (\cite{Demazure; Gabriel 1970}, Chapter II, \S 7, Theorem 3.5),
this actually means  $A[F]\subset\Kernel(\id-h)$.
The Isomorphism Theorem gives an endomorphism $f:A\ra A$   with $\id-h =f\circ F$, or equivalently
$$
h=\id -f\circ F.
$$
Fix a complex embedding $\QQ_\ell^\alg\subset\CC$ and  consider the effect of $f$ and $F$ on the    vector space
$V=H^1(A\otimes k^\alg,\QQ_\ell)\otimes\CC$, which has  dimension $2g$.
Since $f$ and $F$ commute, there is a basis $e_1,\ldots,e_{2g}\in V$ such that the resulting matrices for   $f^*$ and  $F^*$ are
both lower triangular (\cite{A 4-7}, Chapter VII, \S5, No.\ 9, Proposition 19).
Let $\alpha_1,\ldots,\alpha_{2g}$  and $\beta_1,\ldots,\beta_{2g}$ be the     matrix entries on the diagonal for $f^*$ and $F^*$, respectively.
In turn, the matrix for $h^*$ is lower triangular as well, and $1-\alpha_i\beta_i$ are its diagonal entries. These are the
eigenvalues for $h^*$, all of which are   $r$-th roots of unity.
If all   of them  are trivial, then $h^*$ is the identity, because it is diagonalizable, and
hence $h$ is the identity (\cite{Mumford 1970}, Section 18, Theorem 3), contradiction.

Fix an   index $1\leq d\leq 2g$ for which the $r$-th root of unity  $\zeta= 1-\alpha_d\beta_d$   primitive.
Since the characteristic polynomial for $h^*$ belongs to $\ZZ[T]$, all conjugates of $\zeta$
appear among the $1-\alpha_i\beta_i$.
Choose indices $i_1,\ldots,i_{r-1}$ so that
$$
\{1-\alpha_{i_1}\beta_{i_1}, \ldots, 1-\alpha_{i_{r-1}}\beta_{i_{r-1}}\}=\{\zeta^1,\ldots,\zeta^{r-1}\}.
$$
Using this and the cyclotomic polynomial $P(T)=T^{r-1}+\ldots  + 1 = \prod_{j=1}^{r-1} (T-\zeta^j)$   for substitutions, we obtain
$$
r =  P(1) =  \prod_{j=1}^{r-1}(1-\zeta^j) = \prod_{j=1}^{r-1}\alpha_{i_j}\cdot \prod_{j=1}^{r-1}\beta_{i_j}.
$$
Now  recall that the $\beta_{i_j}$ are
Weil numbers, having $|\beta_{i_j}|=p^{1/2}$. Taking absolute values from the above equation gives $r=s p^{(r-1)/2}$,
with the factor $s= \prod_{j=1}^{r-1}|\alpha_{i_j}|$. Each $\alpha=\alpha_{i_j}$ is an algebraic integer, so
the same holds for the conjugate $\bar{\alpha}$ and the absolute value $|\alpha|=(\alpha\cdot\bar{\alpha})^{1/2}$.
Consequently,  the real number $s=rp^{(1-r)/2}$ is also an algebraic integer.

In the special case $r=2$ the factor becomes  $s=2p^{-1/2}$, which is  a root for $T^2-4/p\in\QQ[T]$. 
This polynomial is irreducible regardless of $p$,
hence must be the minimal polynomial for $s$. The latter is an algebraic integer, so $4/p$ must be an integer,
and thus $p=2=r$.
In the general case $r\geq 3$  our algebraic integer $s=rp^{(1-r)/2}$ already belongs to $\QQ$, and thus
is contained in $\ZZ$. From  the uniqueness of prime factorization, we infer $r=p$ and $s=1$
and $(r-1)/2=1$. The latter ensure $p=r=3$. This establishes (i) and (ii).

It remains to verify (iii). Set $N=\Kernel(h-\id)$ and $A'=A/N$.
Recall that by  \cite{Yu 2013}, Theorem 2.9 we have to verify that, with respect to all complex embeddings $\QQ_\ell\subset\CC$,
the eigenvalues of Frobenius on $H^1(A'\otimes k^\alg,\QQ_\ell)\otimes\CC$ take  the form $\xi\cdot p^{1/2}$ for some
root of unity $\xi\in\CC^\times$, and thus are supersingular Weil numbers.

Let us first reduce to the case that $N$ is finite. Write $q:A\to A/N=A'$ 
for the quotient map. Since the commutative group schemes of finite type form an abelian category,
there is a unique  monomorphism $i:A'\to A$ such that $i\circ q=h-\id$. Our
$h\in\Aut(A)$ induces by construction an automorphism $h'\in\Aut(A')$ of finite order, which satisfies the equation  
$q\circ (h-\id)=(h'-\id)\circ q$. We claim that $h'-\id$ has a finite kernel. Suppose this is not the case.
To produce a contradiction, we may replace $k$ by its algebraic closure. Using the surjectivity of $q:A\ra A'$,
we find a closed point $a\in A$ such that $q(a)$ has infinite order and belongs to $\Kernel(h'-\id)$.
Then $(h-\id)(a)\in N$, and thus $(h-\id)^2(a)=0$.
On the other hand, the gcd of $(T-1)^2$ and 
$T^r-1$ in $\QQ[T]$ is $T-1$, so there is an integer $m\geq 1$ and polynomials $Q(T), R(T)\in\ZZ[T]$ such that 
$Q(T)\cdot (T-1)^2+R(T)\cdot (T^r-1)=m(T-1)$.
We conclude   
$$
m\cdot(h-\id)(a)=(h-\id)(m\cdot a)=(i\circ q)(m\cdot a)=i(m\cdot q(a))=0.
$$ 
Since $i$ is a closed immersion, we see that $q(a)\in A$ is a torsion point,   contradiction. 
We may thus replace $A$ (resp.\ $h$) by $A'$ (resp.\ $h'$) and   suppose that $N$ is finite. 

From     $0=(h-\id)\circ (h^{p-1}+h^{p-2}+\ldots+\id)$
we infer that the all eigenvalues  $\zeta_i=1-\alpha_i\beta_i$ for $h^*$ are primitive
$p$-th roots of unity. Now we can exploit $p=3$:
For every complex embedding $\QQ_\ell\subset\CC$ we   have $\zeta_i=e^{\pm 2\pi i/3}$.
Set $\omega= e^{2\pi i/6}$.
Using $- e^{2\pi i/3}=\omega^{-1}$ and $-e^{-2\pi i/3}=\omega$ and $\omega^{-1}+\omega=1$ we compute 
$$
|1-\zeta_i|^2=(1-e^{2\pi i/3})(1-e^{-2\pi i/3}) =3,
$$
and see $1-\zeta_i = \omega^{\pm 1}\cdot p^{1/2}$.
In particular $|1-\zeta_i|=p^{1/2}$, which coincides with $|\beta_i|=p^{1/2}$.
From $\alpha_i\beta_i = 1-\zeta_i $  we get  $|\alpha_i|=1$. 
Since this  applies to all complex embeddings, the $\alpha_i\in\CC$ are algebraic integers all whose conjugates  lie on the unit circle, so
they actually must be roots of unity  (\cite{Washington 1997}, Chapter I, Lemma 1.6). 
Setting $\xi_i= \alpha_i^{-1}\omega^{\pm 1}$, we get $\beta_i=(1-\zeta_i)/\alpha_i=\xi_i\cdot p^{1/2}$, as desired.
 
\medskip\noindent
\textbf{Step 2:} \emph{Suppose that the ground field $k$ is finite.}
Let $k_0=\FF_p$ be the prime field. Then $k_0\subset k$ is a finite Galois extension, with cyclic $G=\Gal(k/k_0)$,
say of order $m\geq 1$. We now form the Weil restriction $A_0=\Res_{k/k_0}(A)$,
which comes with an induced automorphism $h_0=\Res_{k/k_0}(h)$. By Lemma \ref{weil restriction} we have 
\begin{equation}
\label{base-change weil restriction}
A_0\otimes_{k_0} k = \prod_{\sigma\in G} A_\sigma\quadand h_0\otimes_{k_0} k=\prod_{\sigma\in G} h_\sigma,
\end{equation} 
where $A_\sigma$ denotes the abelian variety $A\otimes_kk$, with base-change   via $\sigma:k\ra k$.
This shows that $A_0$ is   para-abelian   of dimension $g_0=mg$.
Using the image $e_0\in A_0$ of the origin $e\in A$, it becomes an abelian variety, with $h_0\in\Aut(A_0)$.
Then $h_0\neq\id$, while the   induced action on $\Lie(A_0)$ is trivial, because this holds after base-change to $k$.
Using step 1 with $h_0\in \Aut(A_0)$, we immediate see $p\leq 3$. Furthermore, the order of $h_0$ must be a $p$-power,
so by \eqref{base-change weil restriction} this also holds for $h$.  Finally, $A'_0=A_0/\Kernel(h_0-\id)$ is supersingular, which
carries over to 
$$
A'_0\otimes_{k_0} k= \prod_{\sigma\in G} A_\sigma/\Kernel(h-\id).
$$
Then each factor, and in particular $A/\Kernel(h-\id)$, is supersingular.

\medskip\noindent
\textbf{Step 3:} \emph{Now the ground field $k$ is general.}
Let $R_\lambda\subset k$, $\lambda \in L$ be the ordered set of all subrings that are finitely generated over the prime field $\FF_p$.
By \cite{EGA IVc}, Theorem 8.8.2, we find some member $R=R_\lambda$  such that $A$ and $h$ arise  from 
a family of abelian varieties $\foA\ra \Spec(R)$  and some relative automorphism $h:\foA\ra\foA$, which we denote
by the same letter. Set $r=\ord(h_\eta)$. Since the generic fiber $\foA_\eta$ is schematically dense, we actually have $h^r=\id_\foA$.
Localizing further, we may assume that $\Kernel(h^i-\id_\foA)$, $0\leq i\leq r-1$ are flat.
This ensures that for all  
$s\in S$ the element $h_s\in\Aut(\foA_s)$ also has order $r$.
Moreover, the quotient $\foB=\foA/\Kernel(h-\id_\foA)$
exists as an algebraic  space (see for example \cite{Laurent; Schroeer 2024}, Lemma 1.1), which  is actually a family of abelian varieties.

The closed points $s\in S$ are Zariski dense, because $S$ is a Jacobson space, 
and the residue fields $\kappa(s)$ must by finite, by Hilbert's Nullstellensatz.
Applying step 2 with $h_s\in\Aut(\foA_s)$, we see that the characteristic satisfies $p\leq 3$, that the common order $r$   is a $p$-power,
and that the closed fiber  $\foB_s$ is supersingular.
In turn, for the finite flat group scheme $\foB[p]$ the closed fibers are geometrically connected.
By \cite{EGA IVc}, Theorem 9.7.7 the generic fiber is geometrically connected as well,
which implies that $\foB_\eta$ is supersingular.
\qed
 
\medskip
We next study automorphisms of the smooth proper scheme $P$    obtained from the abelian variety $A$ by forgetting the group law, and express
various phenomena in terms of group cohomology. First note that 
$$
\Aut(P)= A(k)\rtimes \Aut(A),
$$
where the semidirect product is formed with respect to the canonical action of $\Aut(A)$ on the abelian group $A(k)$.
We actually have  $\Aut_{P/k}=A\rtimes \Aut_{A/k}$ as group schemes.
Note that the group law for $\Aut(P)$ and its action on $P$ are given by 
$$
(a,h)\cdot (a',h')= (a+h(a'),hh')\quadand (a,h)\cdot x=  a+ h(x).
$$
We now fix some element $h\in\Aut(A)$ of finite order. Each $a\in A(k)$ gives an   element $g=(a,h)$ from $\Aut(P)$.
Then $r=\ord(h)$ divides $\ord(g)$, and one immediately computes
$$
g^i= (a+h(a)+\ldots+h^{i-1}(a), h^i)\quad\quadand\quad g\cdot x=x  \Longleftrightarrow a=x-h(x).
$$
This shows:

\begin{lemma}
\mylabel{order and freeness}
The element $g=(a,h)$ from $\Aut(P)$ has order $r=\ord(h)$ if and only if   $a\in\Kernel(\id+h+\ldots+h^{r-1})$.
In this situation, the scheme of fixed points for  $g:P\ra P$ is empty  if and only if $a\not\in\Image(\id-h)$.
\end{lemma}

More precisely, $g:P\ra P$ has no rational fixed point if and only if $a$ is not the image of $\id-h:A(k)\ra A(k)$.

Consider the cyclic group $G=\ZZ/r\ZZ$, and recall that for each $G$-module $M$ the first group cohomology can be expressed
as 
$$
H^1(G,M)=\Kernel(1+\sigma+\ldots+\sigma^{r-1})/\Image (1-\sigma),
$$
where $\sigma:M\ra M$ is the effect of the canonical generator in $G$.
Using the  $G$-module $M=A(k)$, we rephrase and refine the  above lemma as follows:

\begin{proposition}
\mylabel{cohomological formulation}
The map $(a,h)\mapsto [a]$ gives a bijection between the set of $A(k)$-conjugacy classes in 
\begin{equation}
\label{conjugacy classes}
\{g\in\Aut(P)\mid \text{$\ord(g)=r$, and $g=(a,h)$ for some $a\in A(k)$}\},
\end{equation} 
and the cohomology group $ H^1(G,A(k))$. Moreover, the automorphism $g:P\ra P$ has no rational fixed points
if and only if  $[a]\neq 0$.
\end{proposition}

\proof
By Lemma \ref{order and freeness}, for each element $g=(a,h)$ from the set \eqref{conjugacy classes} the entry $a$ belongs to the kernel
of $1+h+\ldots+h^{r-1}$, and thus gives a cohomology class $[a]\in H^1(G, A(k))$. Moreover,  each cohomology class arises in this way.
The conjugacy action in question is given by
$(b,\id)\cdot(a,h)\cdot(b,\id)^{-1} = (b + a -h(b), h)$,
which shows that the map $(a,h)\mapsto [a]$ becomes injective when passing to  $A(k)$-conjugacy classes.
The statement on the fixed points also follows from the lemma.
\qed

\begin{corollary}
\mylabel{conjugate automorphisms}
Suppose $k$ is algebraically closed. Then each  $g=(a,h) $ of order $r=\ord(h)$  is conjugate in $\Aut(P)$ to  some
$g'=(a',h)$ where  $a'\in A(k)$   is annihilated by an $r$-power.
\end{corollary}
 
\proof
Set $M=A(k)$, which is a divisible group. Let $M'$ be the torsion submodule, and $M''=M/M'$.
This gives a short exact sequence $0\ra M'\ra M\ra M''\ra 0$ of $G$-modules, and
thus an exact sequence
$$
H^1(G,M')\lra H^1(G,M)\lra H^1(G,M'').
$$
The torsion-free divisible group $M''$ carries a unique structure of a vector space over $\QQ$,
necessarily respected by the $G$-action. Thus $H^1(G,M)$ vanishes, because it is both a $\QQ$-vector space and annihilated
by $|G|=r$. Let $M'_0\subset M'$ be the subgroup of elements annihilated by some $r$-power, and $M'_1=M'/M'_0$.
Arguing as above we see that $H^1(G,M'_0)\ra H^1(G,M')$ is surjective, and the statement follows from the proposition.
\qed

\medskip
The following assertion will be useful as well:

\begin{proposition}
\mylabel{addition map}
The addition  map of commutative group schemes
$$
\Kernel(\id+h+\ldots+h^{r-1})\times \Kernel(\id-h)\lra A
$$
is surjective, and its kernel is isomorphic to a subgroup scheme inside $A[r]$.
\end{proposition}

\proof
Set $N=\Kernel(\id+h+\ldots+h^{r-1})$ and $N'= \Kernel(\id-h)$. The kernel of the the addition map
is isomorphic to $N\cap N'$. This intersection is $h$-stable, and the induced map satisfies  $h=\id$ and $\id+h+\ldots+h^{r-1}=0$,
hence $r\cdot \id=0$, and thus $N\cap N'\subset A[r]$.

For the first assertion it  suffices to verify   $\dim(A)\leq \dim(N)+\dim(N')$,
in light of the preceding paragraph.
Set $B=\Image(\id+h+\ldots+h^{r-1})$.
From the relation $(\id+h+\ldots+h^{r-1})(\id-h) =0$ we see   $B\subset N'$.
This gives $\dim(A) -\dim(N) = \dim(B)\leq \dim(N')$, as desired.
\qed

\medskip
The following consequence will be important in the next section:

\begin{proposition}
\mylabel{non-simple and disconnected}
Suppose $g=(a,h)$   has order $r=\ord(h)$, and that the fixed scheme for $g:P\ra P$ is empty.
Then the abelian variety $A$ is not simple, and the finite group scheme $A[r]$ is disconnected.
\end{proposition}

\proof
We have $r\geq 2$ because $g$ has no fixed points. Thus $h\neq \id$, and furthermore the abelian variety $A'=\Image(\id-h)$
does not contain $a$ by Lemma \ref{order and freeness}. Thus $0\subsetneqq A'\subsetneqq A$, hence  $A$ is not simple.

Seeking a contradiction, we assume that $A[r]$ is connected. It is then a local Artin scheme with residue field $k$,
so we may also assume that  $k$ is algebraically closed.
Being successive extensions, the kernels $A[r^\nu]$, $\nu\geq 0$ are connected as well.
By Proposition \ref{conjugate automorphisms} we may assume that some of these kernels contain $a$. 
Thus $a=0$, hence $g=(0,h)$ fixes the origin $0\in P$, contradiction.
\qed

%===========================================================
\section{\texorpdfstring{$T$}{T}-trivial varieties}
\mylabel{T-trivial}

Let $k$ be a ground field of characteristic $p\geq 0$.
Recall that a  locally free sheaf $\shE$ of rank $n\geq 0$ on some scheme $X$ is called \emph{free}
if it is isomorphic to $\bigoplus_{i=1}^ n\O_X $. The central topic of this paper
are smooth schemes  with free tangent sheaf. In the literature, one often finds the locution ``trivial tangent bundle''.  We find the following terminology useful:

\begin{definition}
\mylabel{t-trivial}
A \emph{$T$-trivial variety} is a smooth proper scheme $X$ with $h^0(\O_X)=1$ such that the tangent sheaf $\Theta_{X/k}=\uHom(\Omega^1_{X/k},\O_X)$ is free.
\end{definition}

Equivalently, the sheaf of K\"ahler differentials $\Omega^1_{X/k}$ is free.
Note that  $T$-trivial varieties are very special
cases of \emph{$K$-trivial varieties}, where $\omega_X=\det(\Omega^1_{X/k})$ is isomorphic to the structure sheaf.
Abelian varieties $A$ have $\Theta_{A/k}=\Lie(A)\otimes_k\O_A$, and are thus examples of $T$-trivial varieties. 
Let us collect some basic property:

\begin{proposition}
\mylabel{properties t-trivial}
\begin{enumerate}
\item If $X$ and $Y$ are $T$-trivial, so is their product $X\times Y$.
\item If $f:X\ra Y$ is a finite \'etale covering of a $T$-trivial variety $Y$, then every irreducible component $X'\subset X$
is a $T$-trivial variety over the field of constants $k'=\Gamma(X',\O_{X'})$.
\item Let $k\subset k'$ be a field extension. Then $X$ is $T$-trivial if and only if the base-change $X'=X\otimes k'$ is $T$-trivial.
\item The Weil restriction $X_0=\Res_{k/k_0}(X)$ of a $T$-trivial variety $X$ over $k$ along a finite separable extension  
$k_0\subset k$ is a $T$-trivial variety over $k_0$.
\end{enumerate}
\end{proposition}

\proof
The first statement follows from $\Omega^1_{(X\times Y)/k}=\pr_1^*(\Omega^1_{X/k})\otimes \pr_2^*(\Omega^1_{Y/k})$.
For (ii) it suffices to treat the case that $X$ is irreducible. Since $X$ is smooth over $k$,
the finite   extension $k\subset k'$ must be separable.
The structure maps $X\ra \Spec(k')\ra\Spec(k)$ induces an exact sequence
$\Omega^1_{k'/k}\otimes_k\O_X\ra \Omega^1_{X/k}\ra \Omega^1_{X/k'}\ra 0$,
where the term on the left vanishes. 
In the exact sequence 
$f^*(\Omega^1_{Y/k})\ra \Omega^1_{X/k}\ra\Omega^1_{X/Y}\ra 0$, the term on the right vanishes, and the other  terms
are  locally free of the same rank. Thus the arrow on the left is  bijective, so $\Omega^1_{X/k'}$ is free,
and  (ii) follows.

Statement (iii) follows from Lemma \ref{descent freeness}.
It remains to check (iv). In light of (iii) it suffices to verify that $X_0\otimes_{k_0}k$ is $T$-trivial.
This indeed follows from Lemma \ref{weil restriction}, together with (i).
\qed

\medskip
As for abelian varieties, morphisms between $T$-trivial varieties have remarkable rigidity properties.
Let $h:X\ra Z$ be a surjective morphism between $T$-trivial varieties, with Stein factorization $Y=\Spec h_*(\O_X)$.
The resulting morphisms are denoted by
$$
\begin{tikzcd} 
X\ar[r, "g"']\ar[rr,bend left=40, "h"] 	& Y\ar[r,"f"']	& Z.\\ 
\end{tikzcd}
$$

\begin{proposition}
\mylabel{morphisms t-trivial}
Suppose in the above setting that the function field extension $k(Z)\subset k(X)$ is separable. Then the following holds:
\begin{enumerate}
\item The morphisms $h$  and $g$ are smooth, and   $f$ is   \'etale.
\item The scheme $Y$ and all fibers $g^{-1}(y)$, $y\in Y$ are $T$-trivial varieties.
\item If $X$ and $Z$ have the same dimension, then $h:X\ra Z$ is \'etale.
\item The canonical  sequence
$0\ra h^*(\Omega^1_{Z/k})\ra \Omega^1_{X/k}\ra\Omega^1_{X/Y}\ra 0$ is split exact, and all terms are free.  
\end{enumerate}
\end{proposition}

\proof
Set $m=\dim(X)$ and $n=\dim(Z)$, and consider the exact sequence
$$
h^*(\Omega^1_{Z/k})\lra \Omega^1_{X/k}\lra \Omega^1_{X/Z}\lra 0.
$$
By assumption, the terms on the left are free, and the map on the left is injective at the generic point of $X$.
According to  Lemma  \ref{categorical view}, the term on the right is free, the map on the left is injective, and the short
exact sequence splits, which already gives (iv). Furthermore, $\rank(\Omega^1_{X/Z})=m-n$. 
It then follows from  \cite{Hartshorne 1977}, Chapter III, Proposition 10.4 that $h:X\ra Z$ is smooth.
If $X$ and $Z$ have the same dimension, the generic fiber is zero-dimensional, hence $h$ is finite, and thus \'etale,
giving (iii).

% We next verify that $h:X\ra Z$ is smooth, and in particular flat.
% Since $X$ and $Y$ are irreducible, the generic fiber $h^{-1}(\eta)$ has dimension $m-n$. By Chevalley's Semicontinuity Theorem (\cite{EGA IVc}, Theorem 13.1.3),
% for every $z\in Z$ the irreducible components of $h^{-1}(z)$ have dimension at least $ m-n$.
% It follows that the fibers $h^{-1}(z)$ are smooth (\cite{EGA IVd}, Proposition 17.15.5), and equi-dimensional 
% of dimension $m-n$. So the morphism $h:X\ra Z$ is equi-dimensional  in the sense
% of \cite{EGA IVc}, Definition 13.2.2, and thus  must be smooth (\cite{EGA IVd}, Corollary 17.5.6).

The smoothness of $h$ ensures  that the Stein factorization $f:Y\ra Z$ is \'etale (\cite{EGA IIIb}, Remark 7.8.10).
So in the exact sequence $ f^*(\Omega^1_{Z/k})\ra \Omega^1_{Y/k}\ra \Omega^1_{Y/Z}\ra 0$, the term on the right vanishes,
and the map on the left is injective (\cite{EGA IVd}, Theorem 17.11.1).
Thus $\Omega^1_{Y/k}$ is free, and it follows that $Y$ is a $T$-trivial variety.
Using that $\Omega^1_{X/Y}$ is free, we infer that the $g^{-1}(y)$ are $T$-trivial varieties, which gives (ii).
Applying the preceding paragraph with $g:X\ra Y$ instead of $h:X\ra Z$,
we see that the former is smooth, which establishes (i). 
\qed

\medskip
A  \emph{para-abelian variety} is a  smooth proper scheme $P$ such that for some field extension $k\subset k'$, 
the base-change $P'=P\otimes k'$ admits the structure of an abelian variety. Such $P$ are $T$-trivial
according to Proposition \ref{properties t-trivial}.
As explained in \cite{Laurent; Schroeer 2024}, Section 5
the subgroup scheme $G\subset \Aut_{P/k}$ that acts trivial on $\Pic^\tau_{P/k}$ is an abelian variety,
and its action on $P$ is free and transitive. Moreover, we have an identification 
$\Pic^\tau_{P/k}=\Pic^\tau_{G/k}$.
The \emph{Serre--Lang Theorem} on abelian varieties (\cite{Mumford 1970}, Section 18) takes the following form:

\begin{proposition}
\mylabel{serre-lang}
Let $P$ be a para-abelian variety, and $Q\ra P$ be a finite \'etale covering with connected total space.
Then $Q$ is a para-abelian variety over the field $k'=H^0(Q,\O_Q)$.
\end{proposition}

For every proper scheme $X$ with $h^0(\O_X)=1$, there is   morphisms $f:X\ra P$ to a para-abelian variety $P$ that is universal
for morphisms to para-abelian varieties (\cite{Laurent; Schroeer 2024}, Corollary 10.5).  
Note that  $\Pic^\tau_{P/k}$ coincides with the maximal abelian subvariety in $\Pic^\tau_{X/k}$.
One calls $P=\Alb_{X/k}$ the \emph{Albanese variety}, and $f$ the \emph{Albanese map}.
Its formation functorial, stable under base-change, and equivariant for the action of the group scheme $\Aut_{X/k}$.
The following  observation is due to Mehta and Srinivas (\cite{Mehta; Srinivas 1987}, Lemma 1.4):

\begin{proposition}
\mylabel{albanese t-trivial}
Let $X$ be a $T$-trivial variety and   $P=\Alb_{X/k}$ be its Albanese variety.
Then the Albanese map $h:X\ra P$ is smooth and  $\O_P\ra h_*(\O_X)$ is bijective. 
\end{proposition}

\proof The  map $H^0(P,\Omega^1_{P/k})\ra H^0(X,\Omega^1_{X/k})$ is injective (\cite{Igusa 1955b}, Lemma 1 or \cite{Mehta; Srinivas 1987}, Lemma 1.3),
hence the Albanese map is surjective, with separable function field extension. So by Proposition 3.3 the map $h:X\ra P$ must be smooth, and its Stein factorization  $f:Q\ra P$ is \'etale. Using  Proposition \ref{serre-lang}
we infer  that $Q$ is para-abelian. From the universal property
of the Albanese variety one sees that $f$ admits a section $s:P\ra Q$.  Thus the field extension $k(P)\subset k(Q)$
is an equality, and it follows that $f:Q\ra P$ is an isomorphism.
\qed

\medskip
The discrepancy between $T$-trivial varieties and para-abelian varieties now can be seen as  a question about Betti numbers:

\begin{proposition}
\mylabel{t-trivial and betti}
For each $T$-trivial variety $X$ we have $ b_1(X)\leq 2 \dim(X)$. Moreover, equality holds if and only if $X$ is para-abelian.  
\end{proposition}

\proof
The Albanese map $h:X\ra P$ is smooth, according to the proposition,  hence $\dim(X)\geq \dim(P)$.
The group scheme $\Pic^0_{X/k}$ is proper (\cite{Laurent; Schroeer 2024}, Proposition 2.3), and thus $\dim(P)=\dim(\Pic_{X/k})$.
Finally, we have  $2\dim(\Pic_{X/k}) = b_1(X)$ by Lemma \ref{betti and picard}. Thus 
$$
b_1(X) = 2\dim(P) \leq  2\dim(X).
$$
The outer terms are equal if and only if $X$ and   $P$ have the same dimension. In this case,
the Albanese map $h:X\ra P$ is \'etale, by   Proposition \ref{morphisms t-trivial}, and thus an isomorphism, in light of Proposition \ref{albanese t-trivial}.
\qed

\begin{corollary}
\mylabel{t-trivial char zero}
In characteristic zero, every $T$-trivial variety is para-abelian.
\end{corollary}

\proof
Hodge theory gives  $H^1(X(\CC),\ZZ)\otimes\CC=H^{1,0}(X)\oplus H^{0,1}(X)$,
where the summands on the right have the same dimension. Using  $h^0(\Omega^1_{X/k})=\dim(X)$ one gets
$b_1(X)=2\dim(X)$, and the assertion follows from the Proposition.
\qed

\medskip
Let $X$ be an $n$-dimensional  $K$-trivial variety. Serre duality ensures $h^n(\O_X)=1$,
and the relative Frobenius map $F:X\ra X^{(p)}$ yields a linear map
$$
F^*:H^n( X^{(p)}, \O_{X^{(p)}}) \lra H^n( X, \O_X) 
$$
between one-dimensional vector spaces. Choosing a non-zero vector $a\in H^n(X,\O_X)$
yields via  $F^*(a^{(p)}) =\lambda a$ some scalar $\lambda\in k$, whose class modulo $k^{\times p-1}$ 
does not depend on the vector. If $\lambda\neq0$, it  is customary to call $X$ \emph{ordinary}.
Note that there are several other, more refined versions, involving   sheaves of cocycles and coboundaries
in the de Rham complex $\Omega^\bullet_{X/k}$. For $T$-trivial varieties, however, all these notions coincide,
and are also equivalent to the condition that $X$ is \emph{Frobenius split}, that is,  $\O_{X^{(p)}}\ra F_*(\O_X)$ admits
a retraction  (\cite{Mehta; Srinivas 1987}, Lemma 1.1).

Let $X$ be a $T$-trivial variety. Mehta and Srinivas suggested that  the collection of all  finite \'etale
coverings $X_\lambda\ra X$, $\lambda\in L$ with connected total space should contain abelian varieties, at least if $k$ is algebraically closed
(\cite{Mehta; Srinivas 1987}, page 191).
Our second main result   pertains to this, 
where we also take into account the possible fields of constants $k_\lambda=H^0(X_\lambda,\O_{X_\lambda})$ over general ground fields:

\begin{theorem}
\mylabel{t-trivial with etale covering}
Let $X$ be a $T$-trivial variety that is not para-abelian, but has  a finite 
surjection $f:P\ra X$  where $P$ is para-abelian over its fields of constants $k'=H^0(P,\O_P)$, and $k(X)\subset k(P)$ is separable.
% *** \'etale eingefügt. S: nicht noetig, das k(P) ueber k(X) separable.
Then the following holds:
\begin{enumerate}
\item The characteristic satisfies $p\leq 3$.
\item The  abelian variety $A=\Aut^0_{P/k'}$ is not geometrically simple and the finite group scheme $A[p]$ is disconnected.
\item For $p=3$ the $T$-trivial variety $X$ is not ordinary, and the abelian variety $A$ has a supersingular quotient.   
\end{enumerate}
\end{theorem}

\proof
The base-change $X'=X\otimes_kk'$ is a $T$-trivial variety over $k'$, and the induced  morphism $f':P\ra X'$ is finite and surjective,
and respects the $k'$-structures.
By Lemma \ref{removing constant extensions}, we may replace $X$ by $X\otimes k'$ and $k$ by $k'$, and 
thus may assume $P$ is a para-abelian variety over $k$.
In light of Lemma \ref{supersingular quotient}, it now  suffices to treat the case that $k$ is algebraically closed. 

Choose a finite \'etale covering $Q\ra P$ with connected total space such that the composite map $Q\ra X$ is Galois.
Then   $Q$ is para-abelian, by the Serre--Lang Theorem.
Fix a rational point $e_Q\in Q$, and write $e_P\in P$ for the image point. One  may regard the pairs
$(Q,e_Q)$ and $(P,e_P)$ as abelian varieties, and the morphism $Q\ra P$ as an isogeny.
According to Lemma \ref{supersingular quotient}, we may replace $P$ by $Q$, 
and thus may assume that $P\ra X$ is Galois.
Write $G=\Aut(P/X)$ for the Galois group, such that $X=P/G$.
Using Grothendieck's spectral sequence for 
equivariant cohomology (\cite{Grothendieck 1957}, Theorem 5.2.1), we see that the canonical inclusion 
$H^0(X,\Omega^1_{X/k})\subset  H^0(P,\Omega^1_{P/k})^G $
is an equality. From   $h^0(\Omega^1_{X/k})=h^0(\Omega^1_{P/k})$ we
infer that  the induced $G$-action on  $H^0(P,\Omega^1_{P/k})$ is trivial.

Fix a rational point $e\in P$. The resulting map $A\ra P$ is an isomorphism, so one may regard $A$ as the pair $(P,e)$,
and obtains  
$\Aut(P)=A(k)\rtimes\Aut(A)$. The translational part $N=G\cap A(k)$ is normal in $G$,
and the quotient $P/N$ remains para-abelian. So without restriction, we may assume
that the projection $G\ra\Aut(A)$ is injective. Note that the group $G$ is non-trivial, because 
$X$ is not para-abelian.
 
Fix a non-trivial element   $g\in G$, and write  $g=(a,h)$ with   $a\in A(k)$ and $h\in \Aut(A)$.
Let $r\geq 2$ be the common order for $g$ and $h$. 
Since  $g$ and the  translation $a$ act  trivially
on $H^0(A,\Omega^1_{A/k})$, the same holds for $h$. In turn, it also acts trivially 
on the dual vector space $H^0(A,\Theta_{A/k})$, and thus  belongs to the kernel
of the canonical map $\Aut(A)\ra \GL(\Lie(A))$.  
By Theorem \ref{p-group}, we must have $p\leq 3$, and $r=\ord(h)$ is some $p$-power.
Using Proposition \ref{non-simple and disconnected}, we see that $A$ is not simple and  contains a point of order $p$.
For $p=3$ we also have a supersingular quotient, again by Theorem \ref{p-group}.
\qed

\medskip
Let us recall Igusa's construction \cite{Igusa 1955} of $T$-trivial surfaces, and observe that the key features carry over to higher dimensions:
Let $E$ and $E'$ be   elliptic curves, $h\in\Aut(E)$ and $a\in E'(k)$ be group elements
both of order $p$. The resulting diagonal action of the cyclic group $G=C_p$ on 
the abelian surface $A=E\times E'$ is free. It turns out that the induced action on $H^0(A,\Omega^1_{X/k})$ is trivial,
and   the quotient $X=A/G$ must be  $T$-trivial. Note that this is a \emph{bielliptic surfaces},
which have   Betti numbers are $b_1= b_2=2$. 
Also note that for  $k=k^\alg$, the   possibilities for the group $\Aut(E)$ are
$$
C_2,\quad C_4,\quad  C_6,\quad  C_3\rtimes C_4,\quadand Q\rtimes C_3,
$$
formed with cyclic groups  $C_i$ and the quaternion group $Q=\{\pm 1,\pm i,\pm j,\pm k\}$.
We see that  Igusa's  construction is  only possible in  characteristic $p\leq 3$, requires
an abelian variety $A=E\times E'$ that is not simple and contains a point of order $p$,
and for $p=3$ has a  supersingular quotient. 

We also record the  following consequence, which  was already established by Li   (\cite{Li 2010}, Theorem 0.3):

\begin{corollary}
\mylabel{li theorem}
Every ordinary $T$-trivial variety $X$   in characteristic $p\geq 3$ is para-abelian.
\end{corollary}

\proof
We may assume that $k$ is algebraically closed. According to the results of Mehta and Srinivas 
(\cite{Mehta; Srinivas 1987}, Theorem 1),
there is a   finite \'etale covering $A\ra X$ by some  abelian variety $A$, which must be ordinary
(\cite{Mehta; Srinivas 1987}, Theorem 1 and Lemma 1.2).
Seeking a contradiction, we assume that $X$ is not para-abelian. The theorem ensure the characteristic is $p=3$,
and the abelian variety $A$ fails to be ordinary, contradiction.
\qed

%===========================================================
\section{Relative Albanese maps}
\mylabel{Relative albanese}

In this section we investigate the existence 
of relative Albanese maps, and what   to do if they fail to exist. The results, which appears to be of independent interest,  will be used in the next
section to reduce the study of general $T$-trivial varieties to those with Betti number $b_1=0$.
Throughout this section, we work over a  fixed base scheme $S$.

The theory of relative Albanese maps was developed in 
\cite{FGA VI},  \cite{Brochard 2021}, \cite{Laurent; Schroeer 2024} and  \cite{Schroeer 2024}.
Let us recall some basic facts.
Suppose  $f:X\ra S$ is flat  proper morphism of finite presentation, and
with $\O_S=f_*(\O_X)$. Then the abelian sheaf $R^1f_*(\GG_m)$, formed with respect to the fppf topology,
is representable by an algebraic space $\Pic_{X/S}$, and the abelian subsheaf given by 
the numerically trivial invertible sheaves is representable by an open subspace $\Pic^\tau_{X/S}$,
which is of finite presentation over $S$.  Note that the formation of these algebraic spaces commutes
with arbitrary base-change. 

A \emph{family of para-abelian variety} is a proper morphism  $P\ra S$ of finite presentation  all whose fibers are para-abelian,
that is, admit the structure of an abelian variety after some ground field extension.
It then follows that the subgroup space $G\subset\Aut_{P/S}$ that fixes $\Pic^\tau_{P/S}$ is a family of abelian varieties,
that its action on $P$ is free and transitive, and that $\Pic^\tau_{P/S}$ is dual to $G$. 
  
A \emph{relative Albanese map } is a morphism $g:X\ra P$ to a family of para-abelian varieties that 
is universal for arrows to families of para-abelian varieties.
Equivalently, the homomorphism $\Pic^\tau_{P/S}\ra\Pic^\tau_{X/S}$ induces, for each $s\in S$, an 
identifications of the abelian variety  $\Pic^\tau_{P/S}\otimes\kappa(s)$
with the \emph{maximal abelian subvariety} $\Pic^\alpha_{X_s/\kappa(s)}$ inside $\Pic^\tau_{X_s/\kappa(s)}$,
compare  \cite{Laurent; Schroeer 2024}, Definition 8.1 and Theorem 10.2. If these conditions hold, we also write $P=\Alb_{X/S}$.
Note that despite uniqueness, the existence of relative Albanese maps in specific situations is often unclear.

If $S$ is integral, with generic point $\eta\in S$, we write
$A_\eta\subset \Pic^\tau_{X_\eta/\kappa(\eta)}$ for  the maximal abelian subvariety,
and    $A\subset \Pic^\tau_{X/S}$ for its schematic closure. It is unclear whether this is a family of subgroup schemes,
let alone a family of abelian varieties. We start with some easy observations:

\begin{lemma}
\mylabel{various relative pics}
Suppose $S$ is integral and noetherian, and $f:X\ra S$ is smooth.  Then the  algebraic spaces $\Pic^\tau_{X/S}$ and $A$ are proper, and $A$ is  equi-dimensional over $S$. Moreover, 
for each   $s\in S$  the following holds:
\begin{enumerate}
\item The affinization  $\Pic^\aff_{X_s/\kappa(s)}$ of the group scheme $\Pic^\tau_{X_s/\kappa(s)}$ is finite.
\item The kernel of the affinization map  is the maximal abelian subvariety  $\Pic^\alpha_{X_s/\kappa(s)}$.
\item We have $\Pic^\alpha_{X_s/\kappa(s)}=(A_s)_\red$ as closed subschemes inside  $\Pic^\tau_{X_s/\kappa(s)}$.
\item The fiber $A_s$ is reduced provided that the group scheme $\Pic^\aff_{X_s/\kappa(s)}$ is reduced.
\end{enumerate}
\end{lemma}

\proof
The algebraic space  $\Pic^\tau_{X/S}$ is proper by \cite{Laurent; Schroeer 2024}, Proposition 2.3, and the same holds for the
closed subspace $A$.  For the remaining statements it suffices to treat the case that $S$ is the spectrum of a complete local noetherian ring $R$,
with closed point $s\in S$ and separably closed residue field.

To see (i) and (ii) we consider the  group scheme $G=\Pic^\tau_{X_s/\kappa(s)}$, which is  of finite type over the  residue field $k=R/\maxid_R$.
By \cite{Demazure; Gabriel 1970}, Chapter III, Theorem 8.2 it sits in a short exact sequence
$0\ra N\ra G\ra G^\aff\ra 0$
where $G^\aff=\Spec\Gamma(G,\O_G)$ is the affinization, and the kernel $N$ is \emph{anti-affine}, which means $h^0(\O_N)=1$.
The latter is an extension of an abelian variety $N/H$ by some smooth connected affine group scheme $H$, according to \cite{Brion 2017}, Lemma 3.1.4. 
In our situation, $G$ is proper, whence $G^\aff$ is finite and $H$ is trivial. This gives (i) and (ii).

Since $A$ is irreducible, the closed fiber $A_s$ must be connected, according to  Hensel's Lemma
(\cite{EGA IVd}, Theorem 18.5.11). Since $A$ contains the zero section of $P$, we thus get $(A_s)_\red\subset \Pic_{X_s/\kappa(s)}^\alpha$, by (i) and (ii).
Set $P=\Pic^\tau_{X/S}$,   suppose for the moment that
the function  $s\mapsto\dim(P_s)$ is constant,  and   write  $g\geq 0$ for the common value. With Chevalley's Semicontinuity Theorem (\cite{EGA IVc}, Corollary 13.1.5) we get
$g\geq \dim(A_s)\geq \dim(A_\eta)=g$, so $A$ is equi-dimensional.
Moreover, the inclusion  $(A_s)_\red\subset \Pic_{X_s/\kappa(s)}^\alpha$ is an equality, because both  schemes are $g$-dimensional, and the right hand side is irreducible.
This yields (iii).  To see (iv), suppose   that $P^\aff_s$ is reduced. Then $P_s^\alpha$ is
a connected component of $P_s$, and thus the inclusion $P_s^\alpha\subset A_s$ must be an equality. This gives (iv).

It remains to check that $s\mapsto\dim(P_s)$ is constant, which by Chevalley's Semicontinuity Theorem
boils down to $\dim(P_s)\leq \dim(P_\eta)$.
Fix a prime $\ell$ that does not divide the order of the finite group scheme $P_s^\aff$ and $P_\eta^\aff$,
and also differs  from the characteristic $p\geq 0$ of the residue field $ R/\maxid_R$.
 Consider the infinitesimal neighborhoods $X_n=X\otimes R/\maxid_R^{n+1}$. 
The short exact sequence $0\ra \O_{X_0}\ra \O^\times_{X_{n+1}}\ra\O^\times_{X_n}\ra 1$ induces an exact sequence
$$
H^1(X_0,\O_{X_0})\lra \Pic(X_{n+1})\lra \Pic(X_n)\lra H^2(X_0,\O_{X_0})
$$
Since $\ell$ is prime to the characteristic of the residue field $R/\maxid_R$, the map in the middle induces
a bijection on $\ell$-torsion elements.
It then follows from Grothendieck's Existence Theorem (\cite{EGA IIIa}, Theorem 5.4.1) that the restriction map  $\Pic(X)\ra \Pic(X_s)$ indeed induces a bijection on $\ell$-torsion.
So 
$\dim_{\FF_\ell} \Pic(X_s)[\ell] \leq \dim_{\FF_\ell}\Pic(X_\eta)[\ell]$.
In turn, we have $\dim(P_s)\leq \dim(P_\eta)$.
\qed

\medskip
This already gives  a sufficient condition for the existence of relative Albanese maps:

\begin{theorem}
\mylabel{existence relative albanese}
Suppose $S$ is integral, normal and excellent, $f:X\ra S$ is smooth, and   the fibers for the structure morphism $A\ra S$ are reduced.
Then the relative Albanese map $g:X\ra\Alb_{X/S}$ exists.
\end{theorem}

\proof
We first check that $A\ra S$ is a family of para-abelian varieties. Being the closure of the regular scheme $A_\eta$,
the total space $A$ contains   no embedded components. By Lemma \ref{various relative pics}, the structure morphism
$A\ra S$ is proper and equi-dimensional, and the fiber-wise inclusions   $\Pic^\alpha_{X_s/\kappa(s)}\subset A_s$
are equalities.
According to  Koll\'ar's generalization of Hironaka's Flatness Lemma
(\cite{Kollar 1995}, Corollary 11), the morphism $f:A\ra S$ must be flat, and thus is a family of para-abelian varieties.

Obviously, the zero section $e:S\ra\Pic^\tau_{X/S}$ factors over $A$.
By \cite{Laurent; Schroeer 2024}, Proposition 4.3, there is a unique group law that turns $A$ into a family of abelian varieties.
The inclusion $A\subset\Pic^\tau_{X/S}$ respects the zero section. Using \cite{Mumford; Fogarty; Kirwan 1993}, Corollary 6.4
one infers that it  actually respects the group laws. We already observed that  each $A_s$ is the maximal abelian subvariety
inside $\Pic^\tau_{X_s/\kappa(s)}$, so by  \cite{Laurent; Schroeer 2024}, Theorem 10.2 the relative Albanese map $g:X\ra\Alb_{X/S}$ exists.
\qed

\medskip
This basically settles the case of characteristic zero:

\begin{corollary}
\mylabel{unconditional existence}
Suppose  $S$ is a $\QQ$-scheme that is integral, normal  and excellent, and that $f:X\ra S$ is smooth.
Then the relative Albanese map $g:X\ra\Alb_{X/S}$ exists.
\end{corollary}

\proof
By Cartier's Theorem (\cite{Demazure; Gabriel 1970}, Chapter II, Theorem 1.1),   group schemes of finite type over the points $s\in S$ are  automatically reduced.
Lemma \ref{various relative pics} ensures that the fiber  $A_s$ is reduced, so the theorem applies.
\qed

\medskip
Without the assumption on the characteristic, the conclusions holds true if for all points $s\in S$ we have $H^2(X_s,\O_{X_s})=0$,
because then $\Pic^\tau_{X_s/\kappa(s)}$ and hence also its affinization are smooth (\cite{Mumford 1966}, Lecture 27). 

All the above results, however, are still insufficient for the applications we have in mind. For lack of better existence  criteria,
we seek to weaken the very notion of  relative Albanese maps. Our first observation in this direction is:

\begin{lemma}
\mylabel{extension   pic}
Suppose $S$ is integral and noetherian, and $f:X\ra S$ is smooth. Then 
the abelian variety $A_\eta$ extends to a family of abelian varieties over some open set $U\subset X$
containing all points $s\in S$ where $\O_{S,s}$ is regular of dimension one.
\end{lemma}

\proof
Since $A_\eta$ extends over some dense open subset, it suffices to treat the case that $S$
is the spectrum of a discrete valuation ring $R$, and $s\in S$ is the closed point.
The task is to check that $A_\eta$ has \emph{good reduction}. 
For this we further may assume that $R$ is complete,  with separably closed residue field 
(\cite{Bosch; Luetkebohmert; Raynaud 1990}, Section 7.2, Theorem 1).
 
Fix a prime $\ell>0$ that does not divide the order of the finite group scheme $\Pic^\aff_{X_\eta/\kappa(\eta)}$,
and also differs  from the characteristic $p\geq 0$ of the residue field $ R/\maxid_R$.
The former ensures that  the inclusion $A_\eta[\ell]\subset\Pic^\tau_{X_\eta/\kappa(\eta)}[\ell]$ is an equality.
In light of the 
N\'eron--Ogg--Shafarevich Criterion (\cite{Serre; Tate 1968}, Theorem 1), our task is to verify
that the finite \'etale group scheme  $\Pic^\tau_{X_\eta/\kappa(\eta)}[\ell]$ is constant.
Equivalently, the $\FF_\ell$-vector space $\Pic(X_\eta)[\ell]$ has dimension $2g$, where $g\geq0$
is the relative dimension of $\Pic^\tau_{X/S}$.
It thus suffices to verify that the restriction map $\Pic(X)[\ell]\ra\Pic(X_0)[\ell]$ is surjective,
where $X_0$ is the closed fiber. This indeed follows as in the proof for Lemma \ref{various relative pics}.
\qed

\medskip
According to \cite{Laurent; Schroeer 2024}, Corollary 10.6 there is an open neighborhood $V$ of the generic point $\eta\in S$
such that  $ X_V\ra\Alb_{X_V/V}$ exists. 
A necessary condition for $V=S$ is that
the para-abelian variety $P_\eta=\Alb_{X_\eta/\kappa(\eta)}$ extends
to a family of para-abelian varieties $P$ over  $S$.  
In this situation, the  Albanese map $g_\eta:X_\eta\ra P_\eta$ can be seen  as  rational map $g:X\dashrightarrow P$
between integral noetherian schemes. Its \emph{domain of definition} $\Dom(g)$ is an open set in $X$,
so its image  $U=f(\Dom(g))$  is an open set in $S$. In fact, $U$ comprises  all $s\in S$
where our the rational map   induces a rational map $g_s:X_s\dashrightarrow P_s$ on the fiber. One then  says that $g$ is a \emph{$U$-rational map}.

\begin{proposition}
\mylabel{extension   albanese map}
Suppose that $S$ is integral, normal and noetherian, $f:X\ra S$ is smooth, and  that  the Albanese variety $P_\eta=\Alb_{X_\eta/\kappa(\eta)}$ extends
to a family of para-abelian varieties $P$ over $S$. Then the open set $U=f(\Dom(g))$ contains every 
codimension-one point $s\in S$, and the Albanese map $g_\eta:X_\eta\ra P_\eta$ extends to a morphism $g_U:X_U\ra P_U$.
\end{proposition}

\proof
For the statement on $U$, it suffices to treat the case that $S$ is the spectrum of a discrete valuation ring $R$,
with closed point   $s\in S$.
Our task is to verify that the inclusion $X_\eta\subset \Dom(g)$ is strict.
For this we argue as in \cite{Bosch; Luetkebohmert; Raynaud 1990}, Section 2.5, Proposition 5:
Let $\Gamma\subset X\times P$ be the closure of the graph 
for the Albanese map $g_\eta:X_\eta\ra P_\eta$. We have to check that the projection $\pr_1:\Gamma\ra X$ is an isomorphism.
By fppf  descent, we may replace $S$ by $X$, and assume   that the structure morphism $f:X\ra S$ has a section.
The induced section for $P$   endows it   with  the structure of a family of abelian varieties.  
This  is actually the \emph{N\'eron model} of $P_\eta$. Since $X\ra S$ is smooth,
the Albanese map $g_\eta:X_\eta\ra P_\eta$ extends to a morphism over $S=\Spec(R)$. 

For the remaining statement, we may assume that $\Dom(g)$ surjects onto the normal scheme $S$, such that $g:X\dashrightarrow P$
is an $S$-rational map.
Each codimension-one point $\zeta\in X$  either belongs to the generic fiber $X_\eta$, or maps to a codimension-one point
$s\in S$.
By the previous paragraph,  the rational map $g:X\dashrightarrow P$
is defined at all such $\zeta\in X$.  
By the Weil Extension Theorem (\cite{Bosch; Luetkebohmert; Raynaud 1990}, Section 4.4, Theorem 1), 
the $S$-rational map $g$ is defined everywhere.
\qed

\medskip
Our main result on relative Albanese maps is a weak  form of existence:

\begin{theorem}
\mylabel{weak relative albanese map}
Suppose that  $S$ is integral and normal, $f:X\ra S$ is smooth, and that the generic fiber $X_\eta$ contains a rational point. 
After removing a closed set $Z\subset S$ of codimension at least two,
$P_\eta=\Alb_{X_\eta/\kappa(\eta)}$ extends
to a family of abelian varieties $P$ over $S$, and the Albanese map $ g_\eta:X_\eta\ra P_\eta$ extends to a morphism $g:X\ra P$.
\end{theorem}

\proof
Let $G_\eta=\Pic^\tau_{A_\eta/\kappa(\eta)}$ be the dual for the abelian variety $A_\eta=\Pic^\alpha_{X_\eta/\kappa(\eta)}$.
By Lemma \ref{extension   pic} we may assume that it extends to a family of abelian varieties $G\ra S$.
Fix a rational point on $X_\eta$. The resulting rational point on the Albanese variety yields an identification 
$P_\eta=\Alb_{X_\eta/\kappa(\eta)}$ with $G_\eta$.
Thus $P_\eta$ extends to a family $P\ra S$. By Proposition \ref{extension albanese map}, 
the Albanese map $g_\eta:X_\eta\ra P_\eta$ extends to a morphism $g:X\ra P$,
at least after removing a closed subset $Z\subset S $ of codimension at least two.
\qed

\medskip
Although it is unclear whether the above $g:X\ra P$ enjoys a universal property, it can serve as a useful
substitute for   the relative Albanese map, as we shall see in the next section.

%===========================================================
\section{The first Betti number}
\mylabel{Betti number}

Let $k$ be a ground field of characteristic $p\geq 0$.
Throughout this section,    $X$ denotes  an $n$-dimensional $T$-trivial variety.
To unravel its geometry, the chief tool is the   Albanese variety $B=\Alb_{X/k}$ and
the Albanese map $f:X\ra B$. Recall that the former has $\dim(B)=2b_1(X)$, and that the latter is a family of $T$-trivial varieties,
usually of smaller dimension $n'<n$.
Our third main result  is:

\begin{theorem}
\mylabel{decomposition t-trivial}
Assumptions as above. Then there is a finite \'etale covering $X'\ra X$ whose total space
is a $T$-trivial variety $X'$ over the field $k'=H^0(X',\O_{X'})$ where  the  
fibers of the Albanese map $X'\ra \Alb_{X'/k'}$ are $T$-trivial varieties with 
Betti number $b_1=0$. 
\end{theorem}

In characteristic zero, this is a consequence of the Beauville--Bogomolov Decomposition 
for $K$-trivial varieties (\cite{Bogomolov 1974} and \cite{Beauville 1983}), and then    actually    $X'=\Alb_{X'/k}$.
By the result of Mehta and Srinivas (\cite{Mehta; Srinivas 1987}, Theorem 1), 
this carries over to  positive characteristics provided $X$ is ordinary.
\emph{Our  theorem above raises the question
whether or not $T$-trivial varieties with $b_1=0$ exists, besides the obvious example
of the singleton in dimension $n=0$. And if so, for which primes $p>0$ do they occur?} 
At present, we are unable to offer   further insights on this.
It also would be interesting to know if  one may choose $X'$ without constant field extension.

The proof requires some preparation, and will be  given at the end of the section.
We would like to use  the \emph{relative} Albanese variety $\Alb_{X/B}$ over the \emph{absolute} Albanese variety
$B=\Alb_{X/k}$. As discussed in the previous section, the unconditional existence of such a relative construction is in doubt.
To circumvent this issue, consider commutative diagrams
\begin{equation}
\label{diagram with curves}
\begin{CD}
C'	@>>>	X\\
@VVV		@VVfV\\
C	@>>>	B,
\end{CD}
\end{equation} 
where $C$ and $C'$ are regular curves, the vertical maps are surjective, and
the horizontal maps are finite. Let $\eta\in C$ and $\eta'\in C'$ be the generic points.
For simplicity, we also assume $h^0(\O_{C'})=1$.
By Lemma \ref{extension pic}, the abelian variety $G_\eta$ dual to   the maximal abelian subvariety $A_\eta=\Pic^\alpha_{X_\eta/\kappa(\eta)}$
extends to  a family of abelian varieties $G\ra C$.
Set
$$
G_{\eta'}=G_{\eta}\otimes_{\kappa(\eta)}\kappa(\eta')\quadand G_{C'}=G\times_CC'.
$$
The commutative diagram \eqref{diagram with curves} provides a $\kappa(\eta')$-valued point for $X_\eta$,
which gives an identification $\Alb_{X_{\eta'}/\kappa(\eta')}=G_{\eta'}$.
By Theorem \ref{weak relative albanese map},   the Albanese map   $X_{\eta'}\ra G_{\eta'}$ extends to a morphism
$X_{C'}\ra G_{C'}$.

\begin{proposition}
\mylabel{pic alpha isotrivial}
The family $G\ra C$  of abelian varieties is isotrivial.
\end{proposition}

\proof
It suffices to treat the case that $k$ is separably closed. 
To start with, we also assume that $k$ is algebraically closed.
Consider the invertible sheaf $\omega_{G/C}=e^*(\Omega^g_{G/C})$ on the smooth curve $C$, where  
$e:C\ra G$ denotes the zero section and $g=\dim(G/C)$.
According to \cite{Faltings; Chai 1990}, Chapter V, Proposition 2.2
it suffices to verify   $\deg(\omega_{G/C})\leq 0$.
For this, we may replace $C$ by any finite covering, and assume $C=C'$.

As observed above, the Albanese map $g_\eta:X_\eta\ra G_\eta$ extends to a morphism $g:X\ra G$.
 This yields  an exact sequence
\begin{equation}
\label{kaehler sequence}
g^*(\Omega^1_{G/C})\lra \Omega^1_{X_C/C}\lra \Omega^1_{X_C/G}\lra 0.
\end{equation} 
The term on the left is locally free, and the term in the middle   is free, because $\Omega^1_{X/B}$ is free.
Let $d\geq0 $ be its rank. 
The generic fiber $X_\eta$ is a $T$-trivial variety, hence its Albanese map $g_\eta$ is smooth.
So in \eqref{kaehler sequence} the map on the left is  locally a direct summand over  $X_\eta$,
and in particular injective.
Pulling back along the section $s:C\ra X$ yields an inclusion
$e^*(\Omega^1_{G/C})\subset \O_C^{\oplus d}$, and thus $\deg(\omega_{G/C})\leq 0$.

It remains to cope with the case that $k$ is merely separably closed. Fix a prime $\ell\geq 3$ different from the characteristic $p\geq 0$.
Let $D$ be the normalization of the reduction for $C^\alg=C\otimes k^\alg$.
The canonical map $D\ra C^\alg$ is a finite universal homeomorphism.
Using  \cite{SGA 1}, Expos\'e IX, Theorem 4.10, we find a symplectic level structure $(\ZZ/\ell\ZZ)^{2g}_C\ra G$.
According to \cite{Faltings; Chai 1990}, Corollary 2 for Theorem 6.7, the stack 
$\shA_{g,\ell}$ of $g$-dimensional abelian varieties with such a level structure
is an algebraic space. Let $C\ra \shA_{g,\ell}$ be the classifying map for $G$.
By the previous paragraph, it factors over a singleton after base-change to $k^\alg$, so the same holds over $k$.
\qed

\medskip
Now back to the absolute Albanese variety $B=\Alb_{X/k}$. Fix a prime $\ell>0$ different from
the characteristic $p\geq 0$. Since the Albanese map $f:X\ra B$ is smooth and proper,
the  higher direct images $R^if_*(\QQ_\ell)$ are $\ell$-adic local systems,
of some rank $d_i\geq 0$. 
Taking   fibers over the geometric point $a\in X^\alg$ turns them 
into   representations $\pi_1(B,a)\ra\GL_{d_i}(\QQ_\ell)$,
and the local system is \emph{isotrivial} if and only if the image of the   representation is finite.

\begin{proposition}
\mylabel{first direct image isotrivial}
The $\ell$-adic local system $R^1f_*(\QQ_\ell)$   is isotrivial.
\end{proposition}

\proof
It suffices to treat the case that the ground field $k$ is algebraically closed.  Let $d\geq 0$ be the rank of the local system $R^1f_*(\QQ_\ell)$.
Our task is to show that the  corresponding representation $\pi_1(B,a)\ra\GL_d(\QQ_\ell)$  
has finite image. This is trivial for $n=0$, because then $d=0$, so we  assume $n\geq 1$.
We now introduce particular curves $C$ and $C'$ to form the diagram \eqref{diagram with curves}.

By Bertini's Theorem   together with the Lefschetz Theorem for algebraic fundamental groups 
(\cite{Jouanolou 1983},  Chapter I, Theorem 6.3 and \cite{SGA 2}, Expos\'e XII, Corollary 3.5), 
there is a smooth curve $C\subset B$ containing the image of $a\in X$ such  that the induced map $\pi_1(C,a)\ra\pi_1(B,a)$  is surjective.
Choose an algebraic closure of the function field $\kappa(\eta)=k(C)$. 
Working with the resulting geometric point $\bar{\eta}$ rather than $a$ as base point, we have to show
that the composite map $\rho:\pi_1(C,\bar{\eta})\ra\GL_d(\QQ_\ell)$ has finite image.

By Proposition \ref{pic alpha isotrivial}, the family $g:G\ra C$ is isotrivial. Moreover, the projection $f_C:X_C\ra C$ is smooth.
We thus find a finite branched covering $C'\ra C$ such that
$G_{C'}=G_0\times C'$ for some abelian variety $G_0$ over   $k$, and that $f_{C'}:X_{C'}\ra C'$ has a section.
Choose a section  to get the diagram \eqref{diagram with curves}. Now 
the Albanese map $X_{\eta'}\ra G_{\eta'}$ extends to a morphism $X_{C'}\ra G_{C'}$,
and $g_{C'}:G_{C'}\ra C'$ is given by the projection $G_0\times C'\ra C'$.
  
To restate our findings in terms of fundamental groups, we    lift the geometric point $\bar{\eta} $ to $C'$.
The Proper Base-Change Theorem (\cite{SGA 4c}, Expos\'e XII, Theorem 5.1) gives $R^1(f_C)_*(\QQ_\ell)=R^1f_*(\QQ_l)|C'$, so 
$\rho:\pi_1(C,\bar{\eta})\ra\GL_d(\QQ_\ell)$
is the  representation corresponding to the local system $R^1(f_C)_*(\QQ_\ell)$.  
By Lemma \ref{local systems over ground field}, 
the local systems $R^1(f_{C'})_*(\QQ_\ell)$ and $R^1(g_{C'})_*(\QQ_\ell)$ are isomorphic at the generic point $\eta'$,
and by construction $G\otimes_{\kappa(\eta)}\kappa(\eta')=G_0\otimes_k\kappa(\eta')$.
Thus   the diagram
\begin{equation}
\label{diagram fundamental groups}
\begin{tikzcd}  
\pi_1(\eta',\bar{\eta})\ar[d]\ar[r]	& \pi_1(\eta,\bar{\eta})\ar[r]	& \pi_1(C,\bar{\eta})\ar[d,"\rho"] \\
\pi_1(S,\bar{\eta})\ar[rr,"e"']	&			&  \GL_d(\QQ_l),
\end{tikzcd}
\end{equation} 
is commutative, where $S=\Spec(k)$. Of course, the latter is simply-connected  and the lower horizontal map is trivial.
In the upper row, the map on the left is injective  with image of finite index, and the map on the right  is surjective
by \cite{SGA 1}, Expos\'e V, Proposition 8.2.
It follows that the representation $\rho:\pi_1(B,\bar{\eta})\ra\GL_d(\QQ_\ell)$ vanishes on $\pi_1(\eta',\bar{\eta})$, and thus     has finite image.
\qed

\begin{proposition}
\mylabel{first direct image constant}
Suppose that the $\ell$-adic local system $R^1f_*(\QQ_\ell)$ is constant. Then this   system vanishes, and  the fibers 
of the Albanese map $f:X\ra\Alb_{X/k}$ are $T$-trivial varieties with first Betti number $b_1=0$.
\end{proposition}

\proof
It suffices to treat the case that $k$ is algebraically closed.  
The Leray--Serre spectral sequence for the Albanese map   gives an exact sequence
$$
H^1(B,\QQ_\ell)\lra H^1(X,\QQ_\ell)\lra H^0(B,R^1f_*(\QQ_\ell))\lra H^2(B,\QQ_\ell)\lra H^2(X,\QQ_\ell),
$$
where $B=\Alb_{X/k}$.
By assumption, the term in the middle is a $\QQ_\ell$-vector space of dimension $d=\rank(R^1f_*(\QQ_\ell))$.
The outer maps are injective by Lemma \ref{cohomological injectivity} below. The map on the left is actually bijective,
which follows from   Lemma \ref{betti and picard}. 
Thus $d=0$. According to  the Proper Base Change Theorem,  all geometric fibers of $f:X\ra B$ have    $b_1=0$.
\qed

\medskip
\emph{Proof for Theorem \ref{decomposition t-trivial}.}
Let $X$ be an $n$-dimensional $T$-trivial variety, $B=\Alb_{X/k}$   its Albanese variety, and
$f:X\ra B$ be the Albanese map. 
In light of Lemma \ref{removing constant extensions}, it suffices to treat the case that $k$ is algebraically closed.

The case $n=0$ is trivial, so we   assume $n\geq 1$.
Consider  the collection of all finite \'etale cover $X_\lambda\ra X$, $\lambda\in L$ with connected total space.
Each $X_\lambda$ is an $n$-dimensional $T$-trivial variety, and thus has    $b_1(X_\lambda)\leq 2n$.
Replacing $X$ by some $X_\mu$   whose first Betti number attains the largest value, 
we  may assume   $b_1(X)=b_1(X_\lambda)$ for all $\lambda\in L$.

Set $B=\Alb_{X/k}$, and consider the Albanese map $f:X\ra B$ and the ensuing local system
$R^1f_*(\QQ_\ell)$. The latter is isotrivial, according to Proposition \ref{first direct image isotrivial}.
Choose some finite \'etale cover $B'\ra B$ on which it becomes constant.
Without loss of generality we may assume that $B'$ is connected.
Then $B'$ can be seen as an abelian variety, and   $X'=X\times_BB'$ is an $n$-dimensional $T$-trivial variety.
Let $\pr_2:X'\ra B'$ be the induced projection, and $f':X'\ra\Alb_{X'/k}$
be the Albanese map. Then $\pr_2=g\circ f'$ for some unique  $g:\Alb_{X'/k}\ra B'$.
Since $\pr_2$ and $f'$ are smooth, the same holds for $g$.
By construction  
$$
\dim(\Alb_{X'/k})=2b_1(X')=2b_1(X)=\dim(B)=\dim(B').
$$
Thus the smooth proper morphism $g:\Alb_{X'/k}\ra B'$ is finite.
Since both $\pr_2$ and $g$ are in Stein factorization, the same holds for $g$.
Thus $g$ yields an identification $\Alb_{X'/k}=B'$. By construction, $R^1f'_*(\QQ_\ell)=R^1f_*(\QQ_l)|B'$
is constant.
The assertion now follows from Proposition \ref{first direct image constant}.
\qed

\medskip
The following two observations where used in the proof for Proposition \ref{first direct image constant}:

\begin{lemma}
\mylabel{cohomological injectivity}
Let $Z$ be a smooth proper scheme,   
and $g:Y \ra Z$ be  a proper surjective morphism.   Then the induced maps $H^i(Z^\alg,\QQ_\ell)\ra H^i(Y^\alg,\QQ_\ell)$
are injective, for all $i\geq 0$.
\end{lemma}

\proof
Without loss of generality we may assume $k=k^\alg$ and $h^0(\O_Z)=1$.
Set $n=\dim(Z)$. Choose a closed point $\zeta$ in the generic fiber for $g:Y \ra Z$. By functoriality of cohomology,
we may replace $Y$ be the closure of $\zeta$, and thus may assume that $Y $ is integral and $n$-dimensional.
Let  $\alpha\in H^i(B,\QQ_\ell)$ be a non-zero class. By Poincar\'e Duality there
is a class $\beta\in H^{2n-i}(B,\QQ_\ell)$ such that  $\alpha\cup\beta\neq 0$.
It thus suffices to treat the case $i=2n$. 

Set $\Lambda=\ZZ/\ell^\nu\ZZ$. 
Given an open set $U$ in $Y$, we can consider   cohomology with compact support $H^{2n}_c(U,\Lambda)$.
The inclusion map $i:U\ra Y$ yields a short exact sequence $0\ra i_!(\Lambda_U)\ra \Lambda_Y\ra \Lambda_D\ra 0$,
where $D=Y\smallsetminus U$. In the ensuing long exact sequence
$$
H^{2n-1}(D,\Lambda)\ra H^{2n}(Y,i_!\Lambda)\lra H^{2n}(Y,\Lambda)\lra H^{2n}(D,\Lambda),
$$
the outer terms vanish for dimension reasons, and the second terms coincides with $H^{2n}_c(U,\Lambda)$.
This gives   identification $H^{2n}(Y,\Lambda)=H^{2n}_c(U,\Lambda)$,  compatible with respect to inclusions of open sets.

\newcommand{\Trace}{\operatorname{Tr}}
We now use the trace maps 
$
\Trace_g: R^{2d}g_!\Lambda_U(d)\ra \Lambda_V
$
constructed in \cite{SGA 4c}, Expos\'e XVIII, Theorem 2.9. These satisfy various naturality conditions, and are 
defined for certain morphisms $g:U\ra V$ between separated schemes of finite type.
 The condition is that $g|U_0$ is flat of relative dimension $d$ on some    open set $U_0$,
and   the fibers for $g|U\smallsetminus U_0$  have dimension $<d$.

For some dense open set $V\subset Z$ the preimage $U=f^{-1}(V)$
is smooth and $g=f|U$ is finite and flat. 
This yields a commutative diagram 
 $$ 
\begin{tikzcd}[row sep = tiny]
H^{2n}_c(U,\Lambda(n))\ar[dr]\ar[dd,"\Trace_g"']\\
			& \Lambda\\
H^{2n}_c(V,\Lambda(n))\ar[ur]\\
\end{tikzcd}
$$
of trace maps. Note that by  \cite{SGA 4c}, Expos\'e XVIII,
Theorem 2.14 the diagonal arrows are bijective, so the vertical map is bijective as well.
The composition $\Trace_g \circ g^* $ with the canonical map
$g^*:H^{2n}_c(V,\Lambda(n))\ra H^{2n}_c(U,\Lambda(n))$ is multiplication by $\deg(U/V)$, according to 
loc.\ cit., Theorem 2.9.
Passing to the limit with respect to $\nu$ in  $\Lambda=\ZZ/\ell^\nu\ZZ$ and tensoring with $\QQ_\ell$, 
we infer that $g^*: H^{2n}(Z,\QQ_\ell(n))\ra H^{2n}(Y,\QQ_\ell(n))$
is a bijection of one-dimensional vector spaces over $\QQ_\ell$.
\qed

\begin{lemma}
\mylabel{local systems over ground field} 
Let $Y$ be a geometrically normal proper scheme with $h^0(\O_Y)=1$, and set $P=\Alb_{Y/k}$ and $S=\Spec(k)$.
Write $g:Y\ra S$ and $h:P\ra S$ for the structure maps, and  $f:Y\ra P$ for the Albanese map.
Then    $f^*:R^1 h_*(\QQ_\ell)\ra R^1g_*(\QQ_\ell)$ is a bijection of $\ell$-adic local systems over $S$.
\end{lemma}

\proof
It suffices to treat the case that $k$ is algebraically closed, and  to verify that
$f^*:H^1(P,\ZZ_\ell)\ra H^1(Y,\ZZ_\ell)$ becomes bijective after tensoring with $\QQ_\ell$,
and we may replace the coefficient sheaf $\ZZ_\ell$ by the Tate twist  $\ZZ_\ell(1)$.
Let $A$ be the maximal abelian subvariety in $\Pic^\tau_{Y/k}$.
One may identify $A$ with the maximal abelian subvariety in $\Pic^\tau_{P/k}$, compare
\cite{Laurent; Schroeer 2024}, Definition 8.1 and Proposition 8.3.
Arguing as in the proof for Lemma \ref{betti and picard}, we obtain a commutative diagram
$$
\begin{tikzcd}[row sep = tiny]
		    	& H^1(Y,\ZZ_\ell)\\
\invlim_\nu A[\ell^\nu]\ar[ur]\ar[rd]	\\
		    	& H^1(P,\ZZ_\ell) \ar[uu, "f^*"']
\end{tikzcd}
$$
where the maps from the left to the right are injective, with finite cokernel. Thus $f^*\otimes\QQ_\ell$ is bijective.
\qed

%===========================================================
\section{Liftability}
\mylabel{Liftability}

\newcommand{\MN}{\operatorname{MN}}
Let $k$ be a ground field of characteristic $p\geq 0$, and $X$ be a smooth proper scheme with $h^0(\O_X)=1$,
and of dimension $n=\dim(X)$. 
An important invariant are the \emph{$\ell$-adic Chern classes}
$$
c_i=c_i(X)=c_1(\Omega_{X/k}^1) \in H^{2i}(X,\QQ_\ell(i)),
$$
where $\ell>0$ is a prime different from $p$.
For  $k=\CC$ and $c_1=0$,   the \emph{Beauville--Bogomolov Splitting Theorem}   tells us that there is a finite
\'etale covering $X'\ra X$ such that $X'=A\times Y\times Z$, where $A$ is an abelian variety, $Y$
is  a hyperk\"ahler manifold, and $Z$ is Calabi--Yau variety, at least if $X$ is projective (\cite{Bogomolov 1974} and \cite{Beauville 1983}).
This actually holds for compact K\"ahler manifolds, and the arguments are entirely transcendental. 
% It was recently extended
% to compact Moishezon manifolds by Biswas, Cao, Dumitrecu and Guenancia \cite{Biswas; Cao; Dumitrescu; Guenancia 2024}, and in particular holds without projectivity assumptions 
% in the algebraic setting.
Under the additional assumption $c_2=0$
we  necessarily have $X'=A$, which is already a direct  consequence of Yau's proof of the Calabi Conjecture
(see the discussion in \cite{Kobayashi 1981}). Summing up, $X$ admits a finite \'etale covering by some
abelian variety if and only if $c_1=0$ and $c_2=0$.

From now on we assume   $p>0$. If $X$ admits a finite \'etale covering by some para-abelian variety $P$,
one of course has $c_1=0$ and $c_2=0$. We   now seek to understand to what extend the converse holds true.
Let us say that $X$ \emph{projectively lifts to characteristic zero} if there is a discrete valuation ring $R$ with residue field $R/\maxid_R=k$
and field of fractions $F=\Frac(R)$ of characteristic zero, together with   scheme $\foX$ and 
a projective flat morphism $\foX\ra\Spec(R)$ with closed fiber $\foX\otimes_Rk=X$.

\begin{theorem}
\mylabel{liftings and coverings}
In the above situation, suppose the following holds:
\begin{enumerate}
\item For some $\ell\neq p$, the  $\ell$-adic Chern classes $c_1$ and $c_2$   both vanish.
\item The scheme $X$ projectively lifts to characteristic zero.
\item Characteristic and dimension satisfy  $p\geq 2n+2$.
\item The ground field $k$ is separably closed.
\end{enumerate}
Then there is a finite \'etale covering $A\ra X$ by some  abelian variety $A$.
\end{theorem}

The proof requires some preparation, and will be given at the end of the section. To start with,
we consider for  each integer $n\geq 0$  the expression
$$
\MN(n)=\lim_{t\to\infty}\gcd\Big\{\prod_{i=1}^n(\ell^{2i}-1)\mid\text{$\ell\geq t$ prime}\Big\}.
$$
Note that the integer sequence defining the limit is increasing,  and actually stabilizes.
The eventual value  turns out to be
\begin{equation}
\label{eventual value}
\MN(n)=2^{3n+\val_2(n!)}\cdot\prod_{3\leq \ell\leq 2n+1}\ell^{\nu+\val_\ell(\nu!)},
\end{equation} 
where the product runs over all primes
$\ell$  in the indicated range,   $\nu= \lfloor\frac{2n}{\ell-1}\rfloor$ is a Gau\ss{} bracket,
and $\val_\ell(m)$ denotes the $\ell$-adic valuation. For all this, see \cite{Conrad 2010}, Theorem 6.4 and its proof.
Also note that $\val_\ell(\nu!)= \frac{n-s}{\ell-1}$, where $s = \sum s_i$ is the digit sum in  $\nu=\sum s_i\ell^i$,
see  \cite{Robert 2000}, Chapter 5, Section 3.1.  Let us tabulate the first five values of the above function:
$$
\begin{array}{llllll}
\toprule
n	& 0	& 1		& 2			& 3			& 4\\
\toprule
\MN(n)\quad	& 1	& 2^3\cdot 3\quad 	& 2^7\cdot3^2\cdot5\quad	& 2^{10}\cdot3^3\cdot 5\cdot7\quad	& 2^{15}\cdot 3^4\cdot 5^2\cdot 7\\
\bottomrule
\end{array}
$$
Note   that $\MN(1)=24$ is the order of the group $\{\pm 1,\pm i,\pm j\}\rtimes C_3$, which plays a special role for elliptic curves.
Indeed, our interest in $\MN(n)$ stems from the following general fact, which should be well-known:

\begin{lemma}
\mylabel{bound order}
For each $n$-dimensional abelian variety $A$ and each finite subgroup $G\subset\Aut(A)$, 
the  order $|G|$ divides the integer $\MN(n)$.
\end{lemma} 

\proof 
Fix an ample invertible sheaf $\shL_0$ on $A$.
Then $\shL=\bigotimes_{\sigma\in G}\sigma^*(\shL_0)$ is ample, and its class in $\Pic(A)$ is $G$-fixed.
Write $A^\vee=\Pic^0_{A/k}$ for the dual abelian variety, and consider the homomorphism
$$
f:A\lra A^\vee,\quad a\longmapsto \tau_a^*(\shL)\otimes\shL^\vee,
$$
where $\tau_a(x)=x+a$ denotes translation. By contravariance, the   $G$-action on $A$ induces
an action of the opposite group $G^\op$ on the dual $A^\vee$, which is converted back to  a $G$-action   via 
$\sigma\cdot \shN=(\sigma^{-1})^*(\shN)$.
The above map is  equivariant with respect to these $G$-actions:
First observe
$$
(\sigma^{-1})^*\circ\tau_a^* = (\tau_a\circ\sigma^{-1})^*=(\sigma^{-1}\circ\tau_{\sigma(a)})^*=\tau_{\sigma(a)}^*\circ(\sigma^{-1})^* 
$$
for $\sigma\in G$ and $a\in A(k)$. Together with $(\sigma^{-1})^*(\shL)\simeq\shL$ this gives
$$
\sigma\cdot f(a)=(\sigma^{-1})^*(\tau_a^*\shL)\otimes(\sigma^{-1})^*(\shL^\vee)  =
\tau_{\sigma(a)}^* (\shL) \otimes\shL^\vee=f(\sigma\cdot a).
$$
Furthermore, the homomorphism $f:A\ra A^\vee$ has finite kernel, because the invertible sheaf  $\shL$ is ample. 

Now fix a prime $\ell>0$     not dividing $2p\cdot\deg(f)\cdot|G|$. This has three consequences: First, the induced $G$-equivariant map
$f:A[\ell]\ra A^\vee[\ell]$ is an isomorphism.
With this identification, the  Weil pairing $A[\ell]\times A^\vee[\ell] \ra \mu_\ell$ becomes
a $G$-fixed symplectic form on $A[\ell]$.
Second,  Serre's result ensures that  the homomorphism $G\ra\Aut_{A[\ell]}$ is  injective 
(\cite{Grothendieck 1961}, Appendix or \cite{Mumford 1970}, Section 21, Theorem 5).
Choosing a symplectic basis of $k^\sep$-valued points in the finite \'etale group scheme $A[\ell]$, we obtain an inclusion
$G\subset\Sp_{2n}(\FF_\ell)$. Now recall that the finite   symplectic group   has order
$$
|\Sp_{2n}(\FF_\ell)| = \ell^{n^2}\cdot\prod_{i=1}^n(\ell^{2i}-1),
$$ 
see for example \cite{Grove 2002}, Theorem 3.1. As a third consequence, we see that the order $|G|$ must be a divisor of   $\prod_{i=1}^n(\ell^{2i}-1)$.
This   holds for almost all primes $\ell>0$, and the assertion follows from the very definition of $\MN(n)$.
\qed

\medskip
The relation to $T$-trivial varieties arises as follows:

\begin{proposition}
\mylabel{improvement coverings}
Suppose the ground field $k$ is separably closed.
Let $X$ be a proper scheme of dimension $n\geq 0$, and $f:A\ra X$    be a finite \'etale covering  by some  abelian variety $A$.
Then there is a finite \'etale covering $g:B\ra X$ by another  abelian variety $B$ such that the degree $\deg(B/X)$
divides the integer $\MN(n)$.
\end{proposition}

\proof
Without loss of generality we may assume that $f:A\ra X$ is Galois. Write $P$ for the underlying para-abelian variety
of $A$. The finite group $H=\Aut(P/X)$ is a subgroup of the semi-direct product
$\Aut(P)=A(k)\rtimes\Aut(A)$. Set $N=H\cap A(k)$. The quotient $B=A/N$ is another abelian variety,
and the induced map $g:B\ra X$ is a finite \'etale Galois covering, with relative automorphism group $G=H/N$.
By construction, we have $G\subset\Aut(B)$, and Lemma \ref{bound order} ensures that  $\deg(B/X)=|G|$ divides
the integer $\MN(n)$.
\qed

\medskip
\emph{Proof of Theorem \ref{liftings and coverings}.}  
By assumption, we have a discrete valuation ring $R$ whose residue field $k=R/\maxid_R$ is separably closed of characteristic $p>0$,
and whose field of fraction $F=\Frac(R)$ contains the rational numbers, together with
an algebraic space $\foX$ and a proper flat morphism
to  $S=\Spec(R)$ whose closed fiber $X=\foX\otimes k$ is smooth of dimension $n\geq 0$, with $h^0(\O_X)=1$,
and Chern classes $c_1=0$ and $c_2=0$. Furthermore $p\geq 2n+2$.
Our task is to find a finite \'etale covering $A\ra X$ by some
abelian variety $A$. 

We  may assume that $R$ is henselian, and   contained in the field of complex numbers.
Furthermore, by \cite{SGA 1}, , Expos\'e IV, Theorem 4.10 it suffices to construct the covering after passing to a finite   extension $k'$ of our separably closed
field $k$, and we are thus free to replace $R$ by the normalization in some finite   extension $F\subset F'$.
   
The $\ell$-adic local systems  $R^ih_*(\QQ_\ell)$ are constant, because $R$ is strictly henselian,  and it follows
that the relative Chern classes
$c_i(\Omega^1_{\foX/S}) $ vanishes for $i\leq 2$. So  the generic fiber $U=\foX\otimes_RF$ and
 the complex fiber $V=\foX\otimes_R\CC$ both satisfy $c_1=0$ and $c_2=0$. 
As discussed at the beginning of the section, 
there is a finite \'etale covering $V'\ra V$ by some abelian variety $V'$.
In light of the preceding paragraph and  \cite{SGA 1}, Expos\'e X, Corollary 1.8  we may assume that it arises from $U$ by base-change. 
With Proposition \ref{improvement coverings} we obtain a finite \'etale covering $U'\ra U$
by some abelian variety $U'$ with the additional property that   $\deg(U'/U)$ divides
the integer $\MN(n)$. The latter is relatively prime to $p$, by Formula \eqref{eventual value} and our assumption $p\geq 2n+2$.
Hence, by the theory of specialization of algebraic fundamental group (\cite{SGA 1}, Expos\'e X, Corollary 3.9),
we see that $U'\ra U$ extends to some finite \'etale covering $\foA\ra\foX$.
The schematic closure of the origin $e\in U'(F)$ defines a section $e:S\ra\foA$.
This turns $\shA\ra S$ into a family of abelian varieties, by  \cite{Mumford; Fogarty; Kirwan 1993}, Theorem 6.14,
and the closed fiber yields the desired $A\ra X$.
\qed

\medskip
Note that the projectivity assumptions in Theorem \ref{liftings and coverings} are superflous, and 
the above arguments carry over if $X$ and the structure map $\foX\ra\Spec(R)$ are merely proper, and  $\foX$ is an algebraic space,
such that the complex fiber  $V=\foX\otimes_R\CC$ is a proper algebraic space over $\CC$. Indeed:
By Artin's result (\cite{Artin 1970}, Theorem 7.3), the category of algebraic spaces proper over $\CC$ is equivalent to the category
of compact Moishezon spaces, and the  Beauville--Bogomolov Splitting Theorem was  recently extended
to compact Moishezon manifolds by Biswas, Cao, Dumitrecu and Guenancia \cite{Biswas; Cao; Dumitrescu; Guenancia 2024}.
However, the  finite \'etale covering $A\ra X$ by an abelian variety a posteriori reveals that $X$ must be projective.

%===========================================================


\begin{thebibliography}{ccccc}

\bibitem{Artin 1970}
M.\ Artin:
Algebraization of formal moduli II: Existence of modifications.
Ann.\ Math.\  91 (1970), 88--135.


\bibitem{SGA 4c}
M.\ Artin, A.\ Grothendieck, J.-L. Verdier (eds.):
Th\'eorie des topos et cohomologie \'etale des sch\'emas (SGA 4) Tome 3.
Springer, Berlin, 1973.

\bibitem{Bauer; Gleissner 2022}
I.\ Bauer, C.\ Gleissner: 
Towards a classification of rigid product quotient varieties of Kodaira dimension 0. 
Boll.\ Unione Mat.\ Ital.\ 15 (2022), 17--41. 

\bibitem{Beauville 1983}
A.\ Beauville:
Vari\'et\'es K\"ahleriennes dont la premi\`ere classe de Chern est nulle.
J.\ Differential Geom.\ 18 (1983), 755--782.

\bibitem{Biswas; Cao; Dumitrescu; Guenancia 2024}
I.\ Biswas, J.\  Cao, S.\  Dumitrescu, H.\ Guenancia:
Geometry of $K$-trivial Moishezon manifolds : decomposition theorem and holomorphic geometric structures.
\arXiv{arXiv:2306.16729}, to appear in Math.\ Ann.

\bibitem{Bogomolov 1974}
F.\ Bogomolov:
On the decomposition of K\"ahler manifolds with a trivial canonical class. 
Math.\ USSR Sbornik 22 (1974), 580--583. 

\bibitem{Bombieri; Mumford 1977}
E.\ Bombieri, D.\ Mumford:
Enriques' classification of surfaces in char.\ $p$,  II.
In: 
W.\ Baily, T.\ Shioda (eds.), Complex analysis and algebraic geometry, pp.\ 23--42.
Cambridge University Press, London, 1977.

\bibitem{A 4-7}
N.\ Bourbaki:
Algebra II. Chapters 4--7.
Springer, Berlin, 1990.

\bibitem{Bosch; Luetkebohmert; Raynaud 1990}
S.\ Bosch, W.\ L\"utkebohmert, M.\ Raynaud:
N\'eron models.
Springer, Berlin, 1990.
 
\bibitem{Brion 2017}
M.\ Brion:
Some structure theorems for algebraic groups.
In:
M.\ Can (ed.),  Algebraic groups: structure and actions, pp.\ 53--126.
Amer.\ Math.\ Soc., Providence, RI, 2017. 

\bibitem{Brochard 2021}
S.\ Brochard:
Duality for commutative group stacks.  
Int.\ Math.\ Res.\ Not.\ IMRN 2021,  2321--2388. 

\bibitem{Brown 1982}
K.\ Brown:
Cohomology of groups. 
Springer, Berlin, 1982.

\bibitem{Catanese; Demleitner 2018}
F.\ Catanese, A.\ Demleitner: 
Hyperelliptic threefolds with group $D_4$, the dihedral group of order 8.
Preprint, \arXiv{arXiv:1805.01835}, 2018. 

\bibitem{Catanese; Demleitner 2020}
F.\ Catanese, A.\ Demleitner: 
The classification of hyperelliptic threefolds. 
Groups Geom.\ Dyn.\ 14 (2020),   1447--1454. 

\bibitem{Catanese; Demleitner 2023}
F.\ Catanese, A.\ Demleitner: 
The classification of rigid hyperelliptic fourfolds. 
Ann.\ Mat.\ Pura Appl.\ 202 (2023),  1425--1450. 

\bibitem{Conrad 2010}
B.\ Conrad: Semistable reduction for abelian varieties.
Lecture notes for the Number theory learning seminar 2010--2011.
\url{https://citeseerx.ist.psu.edu/document?repid=rep1&type=pdf&doi=2a1f2325a607872d7a14d70653ef2bb59af20f10}.

\bibitem{Conrad; Gabber; Prasad 2010}
B.\ Conrad, O.\ Gabber, G.\ Prasad:
Pseudo-reductive groups. 
Cambridge University Press, Cambridge, 2010.

\bibitem{Debes; Douai 1997}
P.\ D\'ebes, J.-C.\ Douai:
Algebraic covers: field of moduli versus field of definition. 
Ann.\ Sci. \'Ecole Norm.\ Sup.\ 30 (1997), 303--338.

\bibitem{Deligne 1980}
P.\  Deligne:
La conjecture de Weil. II.
Inst.\ Hautes \'Etudes Sci.\ Publ.\ Math.\  52 (1980), 137--252.

\bibitem{Demazure; Gabriel 1970}
M.\ Demazure, P.\ Gabriel:
Groupes alg\'ebriques.
Masson, Paris, 1970.

\bibitem{Dieudonne 1957}
J.\ Dieudonn\'e:
Groupes de Lie et hyperalgèbres de Lie sur un corps de caractéristique $p>0$. VII.  
Math.\ Ann.\ 134 (1957), 114--133. 

\bibitem{Faltings; Chai 1990}
G.\ Faltings, C.-L.\ Chai:
Degeneration of abelian varieties. 
Springer, Berlin, 1990

\bibitem{Grothendieck 1957}
A.\ Grothendieck:
Sur quelques points d'alg\`ebre homologique. 
Tohoku Math.\ J.\ 9 (1957), 119--221. 

\bibitem{Grothendieck 1961}
% J.-P.\ Serre:
% Rigiditi du foncteur de Jacobi d’echelon n 2 3, Appendix to 
A.\ Grothendieck:
Techniques de construction en geometric analytique. X.
Seminaire Henri Cartan, tome 13,   no.\ 17 (1960/61).

\bibitem{FGA VI}
A.\ Grothendieck:
Technique de descente et th\'eor\`emes d'existence en g\'eom\'etrie alg\'ebrique. VI. 
Les sch\'emas de Picard: propri\'et\'es g\'en\'erales. 
S\'eminaire Bourbaki, Vol.\ 7, Exp.\ 236, 221--243.
Soc.\ Math.\ France, Paris, 1995. 

\bibitem{EGA IIIa}
A.\ Grothendieck:
\'El\'ements de g\'eom\'etrie alg\'ebrique III:
\'Etude cohomologique des faisceaux coh\'erents.
Publ.\ Math., Inst.\ Hautes \'Etud.\ Sci.\ 11 (1961).

\bibitem{EGA IIIb}
A.\ Grothendieck:
\'El\'ements de g\'eom\'etrie alg\'ebrique III:
\'Etude cohomologique des faisceaux coh\'erents.
Publ.\ Math., Inst.\ Hautes \'Etud.\ Sci.\  17 (1963).

\bibitem{EGA IVc}
A.\ Grothendieck:
\'El\'ements de g\'eom\'etrie alg\'ebrique IV: \'Etude locale des
sch\'emas et des morphismes de sch\'emas.
Publ.\ Math., Inst.\ Hautes \'Etud.\ Sci.\  28 (1966).

\bibitem{EGA IVd}
A.\ Grothendieck:
\'El\'ements de g\'eom\'etrie alg\'ebrique IV: \'Etude locale des
sch\'emas et des morphismes de sch\'emas.
Publ.\ Math., Inst.\ Hautes \'Etud.\ Sci.\   32 (1967).

\bibitem{SGA 1}
A.\ Grothendieck:
Rev\^etements \'etales et groupe fondamental (SGA 1).
Soci\'et\'e Math\'ematique de France, Paris, 2003. 

\bibitem{SGA 2}
A.\ Grothendieck:
Cohomologie locale des faisceaux coh\'erents et th\'eor\`emes de Lefschetz locaux et globaux (SGA 2).
North-Holland Publishing Company, Amsterdam, 1968.

\bibitem{Grove 2002}
L.\ Grove:
Classical groups and geometric algebra.
American Mathematical Society, Providence, RI, 2002.

\bibitem{Hartshorne 1977}
R.\ Hartshorne:
Algebraic geometry.
Springer, Berlin,  1977.

\bibitem{Igusa 1955b}
J.-I.\ Igusa:
A fundamental inequality in the theory of Picard varieties.
Proc.\ Nat.\ Acad.\ Sci.\ U.S.A. 41 (1955), 317--320.

\bibitem{Igusa 1955}
J.-I.\ Igusa:
On some problems in abstract algebraic geometry. 
Proc.\ Nat.\ Acad.\ Sci.\ U.S.A. 41 (1955), 964--967.

\bibitem{Ji; Li; McFaddin; Moore; Stevenson 2022}
L.\ Ji, S.\ Li, P.\  McFaddin, D.\  Moore, M.\  Stevenson:
Weil restriction for schemes and beyond. 
In:
P.\ Belman, W.\ Ho, A.\ de Jong (eds.), Stacks Project Expository Collection (SPEC).
Cambridge Univ.\ Press, Cambridge, 2022. 

\bibitem{Joshi 2021}
K.\ Joshi:
On varieties with trivial tangent bundle in characteristic $p>0$.  
Nagoya Math.\ J.\ 242 (2021), 35--51. 

\bibitem{Jouanolou 1983}
J.-P.\ Jouanolou:
Th\'eor\`emes de Bertini et applications.
Prog.\ Math.\  42.
Birkh\"auser, Boston, MA, 1983.

\bibitem{Kobayashi 1981}
S.\ Kobayashi:
Recent results in complex differential geometry.
Jahresber.\ Deutsch.\ Math.-Verein.\ 83 (1981),  147--158. 

\bibitem{Kollar 1995}
J.\ Koll\'ar:
Flatness criteria.
J.\ Algebra 175 (1995),  715--727. 

\bibitem{Lange 2001}
H.\ Lange:
Hyperelliptic varieties.  
Tohoku Math.\ J.\  53  (2001),  491--510.

\bibitem{Laurent; Schroeer 2024}
B.\ Laurent, S.\ Schr\"oer:
Para-abelian varieties and  Albanese maps.
Bull.\ Braz.\ Math.\ Soc.\ 55 (2024), 1--39.

\bibitem{Li 2010}
K.-Z.\ Li:
Differential operators and automorphism schemes. 
Sci.\ China Math.\ 53 (2010),  2363--2380. 

\bibitem{Manin 1963}
Y.\ Manin:
Theory of commutative formal groups over fields of finite characteristic.  
Uspehi Mat.\ Nauk 18 (1963),  (114), 3--90. 

\bibitem{Matsumura 1963}
H.\ Matsumura:
On algebraic groups of birational transformations.
Atti Accad.\ Naz.\ Lincei Rend.\ Cl.\ Sci.\ Fis.\ Mat.\ Nat.\ 34 (1963), 151--155. 

\bibitem{Mehta; Srinivas 1987}
V.\ Mehta, V.\ Srinivas:
Varieties in positive characteristic with trivial tangent bundle.
Compositio Math.\ 64 (1987),  191--212. 

\bibitem{Minkowski 1887}
H.\ Minkowski:
On the theory of positive quadratic forms.  
J.\ für Math.\ CI,  (1887) 196--202. 

\bibitem{Mumford 1966}
D.\ Mumford:
Lectures on curves on an algebraic surface.
Princeton University Press, Princeton, 1966.

\bibitem{Mumford 1970}
D.\ Mumford:
Abelian varieties.
Tata Institute of Fundamental Research Studies in Mathematics 5.
Oxford University Press,  London, 1970.

\bibitem{Mumford; Fogarty; Kirwan 1993}
D.\ Mumford, J.\ Fogarty, F.\ Kirwan:
Geometric invariant theory. Third edition.
Springer, Berlin, 1993.

\bibitem{Novakovic 2012}
S.\ Novakovi\'c:
Absolutely split locally free sheaves on Brauer-Severi varieties of index two.
Bull.\ Sci.\ Math.\ 136 (2012),  413--422. 

\bibitem{Ogus 1979}
A.\ Ogus:
Supersingular $K3$ crystals.
Ast\'erisque 64 (1979), 3--86.

\bibitem{Oort 1974}
F.\ Oort:
Subvarieties of moduli spaces. 
Invent.\ Math.\ 24 (1974), 95--119. 

\bibitem{Pink; Roessler 2004}
R.\ Pink, D.\ Roessler:
A conjecture of Beauville and Catanese revisited. 
Math.\ Ann.\ 330 (2004),  293--308. 

\bibitem{Robert 2000}
A.\ Robert:
A course in $p$-adic analysis. 
Springer, New York, 2000.

\bibitem{Roessler; Schroeer 2022}
D.\ R\"ossler, S.\ Schr\"oer:
Moret-Bailly families and non-liftable schemes.
Algebr.\ Geom.\ 9 (2022),  93--121. 

\bibitem{Schroeer 2024}
S.\ Schr\"oer:
Albanese maps for open algebraic spaces.
Int.\ Math.\ Res.\ Not.\ IMRN 2024,  no.\ 6, 4963--5004. 

\bibitem{Serre; Tate 1968}
J.-P.\ Serre, J.\  Tate:
Good reduction of abelian varieties.
Ann.\ Math.\ \ 88 (1968), 492--517.

\bibitem{Shioda 1978}
T.\ Shioda:
Supersingular $K3$ surfaces.  
In: 
K.\ Lonsted (ed.), Algebraic geometry, pp.\ 564--591.
Springer, Berlin, 1979.

\bibitem{Uchida; Yoshihara 1976}
K.\ Uchida, H.\ Yoshihara:
Discontinuous groups of affine transformations of $\CC^3$.
Tohoku Math.\ J.\  28 (1976), 89--94. 

\bibitem{Washington 1997}
L.\ Washington:
Introduction to cyclotomic fields.  
Springer, New York, 1997.

\bibitem{Yu 2013}
C.-F.\ Yu:
Endomorphism algebras of QM abelian surfaces.  
J.\ Pure Appl.\ Algebra 217 (2013),  907--914. 

\bibitem{Yu 2020}
C.-F.\ Yu:
A note on supersingular Abelian varieties.  
Bull.\ Inst.\ Math.\ Acad.\ Sin.\   15 (2020),  9--32. 

\end{thebibliography}
\end{document}